\documentclass{tac}

\usepackage{amssymb}
\usepackage{amsmath}
\usepackage{array}
\usepackage{bbm}
\usepackage{enumerate}
\usepackage[shortlabels]{enumitem}
\usepackage[latin1]{inputenc}
\usepackage[T1]{fontenc}
\usepackage{shadow}
\usepackage[all,cmtip,2cell]{xy}
\UseAllTwocells
\usepackage{graphicx}
\usepackage{graphbox}
\usepackage{caption}
\usepackage{subcaption}
\usepackage{float}
\usepackage{multicol}
\usepackage[colorlinks=true]{hyperref}
\usepackage{tikz}
\usepackage{tikz-cd}
\usepackage{longtable}
\usepackage{color}
\usetikzlibrary{matrix,arrows,decorations.pathmorphing}
\usepackage{pdfpages}

\newtheorem{conjecture}[thm]{Conjecture}

\newcommand{\CC}{\mathcal{C}}
\newcommand{\QQ}{\mathcal{Q}}
\newcommand{\DD}{\mathcal{D}}
\newcommand{\EE}{\mathcal{E}}
\newcommand{\PP}{\mathcal{P}}
\newcommand{\A}{\mathcal{A}}
\newcommand{\B}{\mathcal{B}}
\newcommand{\W}{\mathcal{W}}
\newcommand{\un}{u}
\newcommand{\co}{c}

\definecolor{myblack}{RGB}{0,0,0}
\definecolor{myred}{RGB}{255,0,0}
\definecolor{myblue}{RGB}{0,0,255}
\definecolor{mygreen}{RGB}{0,120,0}
\definecolor{mypurple}{RGB}{120,0,120}
\definecolor{myorange}{RGB}{255,120,0}
\definecolor{mylightblue}{RGB}{0,255,255}
\definecolor{myolive}{RGB}{120,120,0}
\definecolor{myteal}{RGB}{0,120,120}

\newcommand{\red}{\textcolor{myred}{red}}
\newcommand{\blue}{\textcolor{myblue}{blue}}
\newcommand{\green}{\textcolor{mygreen}{green}}
\newcommand{\purple}{\textcolor{mypurple}{purple}}
\newcommand{\orange}{\textcolor{myorange}{orange}}
\newcommand{\lightblue}{\textcolor{mylightblue}{light blue}}
\newcommand{\olive}{\textcolor{myolive}{olive}}
\newcommand{\teal}{\textcolor{myteal}{teal}}

\DeclareMathOperator{\Id}{id}

\mathrmdef{sk}
\mathrmdef{SW}
\mathrmdef{sw}
\mathrmdef{adj}
\mathrmdef{h}
\mathrmdef{F}
\mathrmdef{Unit}
\mathrmdef{Counit}
\mathrmdef{Sw}
\mathrmdef{Fun}
\mathrmdef{Hom}
\mathrmdef{gSet}
\mathrmdef{id}
\mathrmdef{Map}
\mathrmdef{Adj}
\mathrmdef{Cat}
\mathrmdef{Cusp}

\title{Coherence for adjunctions in a $4$-category}  
\author{Manuel Ara\'{u}jo}     
\address{}
\eaddress{manuel.araujo@tecnico.ulisboa.pt}


\copyrightyear{2024}

\thanks{This work was partially supported by the FCT project grant \textbf{Higher Structures and Applications}, PTDC/MAT-PUR/31089/2017. The work that lead to this paper was also carried out while the author was visiting the \textbf{Max Planck Institute for Mathematics} and later a postdoc with the RTG 1670 \textbf{Mathematics inspired by String Theory and Quantum Field Theory} at the University of Hamburg.}

\begin{document}

\maketitle

\begin{abstract}We give a definition of a coherent adjunction in a $4$-category consisting of a finite list of $k$-morphisms for $k\leq 4$, plus equations beetween $4$-morphisms. We prove that the restriction map from the space of coherent adjunctions in a $4$-category to the space of $1$-morphisms which admit an adjoint is a trivial fibration. We prove that other restriction maps related to fixing parts of the data of an adjunction are also trivial fibrations. We give a conjectural description of a coherent adjunction in an $n$-category.  \end{abstract}

\section{Introduction}

In this paper, we construct a $4$-categorical presentation $\Adj_{(4,1)}$ containing two $0$-cells $X,Y$ and two $1$-cells $l:X\to Y$ and $r:Y\to X$ and we define a \textbf{coherent adjunction} in a strict $4$-category $\CC$ as a functor $\Adj_{(4,1)}\to\CC$. We then prove Theorem \ref{main}, stating that the data of a $1$-morphism $l$ in $\CC$ which admits a right adjoint (by which we mean an adjunction $l\dashv r$ exists in the homotopy $2$-category of $\CC$) can be extended in an essentially unique way to the data of a coherent adjunction.

In order to state this Theorem precisely, denote by $\theta^{(1)}$ the computad consisting of a single $1$-cell, so that $\Map(\theta^{(1)},\CC)$ is the $4$-groupoid of $1$-morphisms in $\CC$, and let $\Map^L(\theta^{(1)},\CC)$ be its full $4$-subgroupoid whose objects are the $1$-morphisms which admit a right adjoint in $\h_2(\CC)$. The map $$E_l:\Map(\Adj_{(4,1)},\CC)\to\Map(\theta^{(1)}, \CC)$$ given by restriction to the $1$-cell $l$ factors through $\Map^L(\theta^{(1)},\CC)$.

\begin{theorem}[Main Theorem]\label{intromain}
	
	Given a strict $4$-category $\CC$, the restriction map $$E_l:\Map(\Adj_{(4,1)},\CC)\to\Map^{L}(\theta^{(1)}, \CC)$$ is a trivial fibration of strict $4$-groupoids.
	
\end{theorem}

The proof builds on an analogous result for adjunctions in $3$-categories proved in \cite{adj3}. Using that result, one is reduced to showing that $$\psi:\Map(\Adj_{(4,1)},\CC)\to\Map(\Adj_{(4,1)},\h_3(\CC))\times_{\Map^L(\theta^{(1)}, \h_3(\CC))}\Map^L(\theta^{(1)}, \CC)$$ is a trivial fibration. This map is a fibration, by a small modification of the main result of \cite{fibrations} (Corollary \ref{corfib}). So we just need to show that its fibres are weakly contractible, which in practice means checking it's surjective on objects and that the homotopy groups of its fibres are trivial. We prove surjectivity by describing how one can lift a coherent adjunction in $\h_3(\CC)$ to one in $\CC$. We prove triviality of homotopy groups by constructing trivialising morphisms for arbitrary elements of these homotopy groups. We do this one cell at a time, by using the lifting properties of fibrations and using the string diagram calculus developed in \cite{data_struct_quasistrict}, \cite{sesquicat}, \cite{sesqui_comp} and \cite{fibrations} for explicit constructions. 

In Section \ref{more_main}, we prove other versions of this Theorem, where instead of fixing only the left adjoint we fix some more of the data constituting a coherent adjunction.

\begin{remark}
	
We restrict to strict $4$-categories as that is the setting where we can rigorously justify the use of our string diagram calculus. In upcoming work, we develop a notion of semistrict $4$-category based on these string diagrams and all results in this paper will then hold in that context, with the same proofs.

\end{remark}

\subsection{A conjectural description of coherent $n$-categorical adjunctions}

The most surprising feature of our description of coherent $4$-categorical adjunctions is that no butterfly relations are needed. This is because given one of the swallowtail $3$-morphisms one can construct the other one, and there is a unique choice that satisfies the butterfly relation (see Section \ref{second_swallowtail}). This means that the swallowtail $3$-cell is the only coherence cell that appears in a minimal description of a coherent $4$-categorical adjunction that is of a genuinely different kind than the ones already appearing in the ususal definition of an adjunction in a $2$-category. This $3$-cell already appeared in the description of a $3$-categorical adjunction. This leads us to conjecture that no more genuinely new coherence is needed in a minimal description of a coherent $n$-categorical adjunction. See Section \ref{conjecture} for more details.

\subsection{Coherence}

We explain in what sense a map $\Adj_{(4,1)}\to\CC$ can be considered a coherent adjunction in $\CC$. This is intented as an informal discussion and we won't use any of it in the rest of the paper. An \textbf{adjunction} between $1$-morphisms in a $4$-category $\CC$ consists of

\begin{description}

 \item[Data:] $1$-morphisms $l,r$ and $2$-morphisms $\un,\co$ in $\CC$.
 \item[Properties:] the two \textit{zigzag} or \textit{snake} composites of $\un$ and $\co$ are equal to appropriate identity $2$-morphisms in the homotopy $2$-category $\h_2(\CC)$.
\end{description}

The properties can be rephrased as follows: \textbf{there exist} $3$-equivalences in $\CC$ from the snake composites of $\un$ and $\co$ to the appropriate identity $2$-morphisms. The fact that these $3$-morphisms are equivalences means that \textbf{there exist} inverse $3$-morphisms and $4$-equivalences from the appropriate composites to the appropriate identity $3$-morphisms. This then means \textbf{there exist} inverse $4$-morphisms such that the appropriate composites are \textbf{equal} to identity $4$-morphisms. 

The fact that we require the \textbf{existence} of certain morphisms as a property, instead of having these morphisms as part of the data of an adjunction is the reason we say this notion of adjunction is \textbf{not coherent}.

Our first attempt at fixing this problem might be to add to the data of an adjunction all the morphisms whose existence we just postulated as part of the properties. This would give a candidate for a definition of a coherent adjunction that would look like

\begin{description}
 \item[Data:] all morphisms specified above.
 \item[Properties:] equations between $4$-morphisms specified above.
\end{description}

This defines a certain coherent $4$-categorical structure, where we use \textbf{coherent} to mean that the structure is defined by certain \textbf{data plus equations at the top dimension}, i.e. a $5$-computad.

However, we still don't know if this new structure is equivalent to that of an adjunction. By this we mean that, given $(l,r,\un,\co)$ in $\CC$ defining an adjunction, there should be an \textbf{essentially unique} way to extend this to the data and equations of our coherent structure. To make this precise, denote by $\PP$ the $5$-computad corresponding to the coherent structure we defined above. Denote by $\Adj_{(2,1)}$ the $3$-computad containing the data ($l$, $r$, $\un$, $\co$) and equations (two snake equations) defining an adjunction in a $2$-category and denote by $\sk_2(\Adj_{(2,1)})$ its $2$-skeleton, obtained by discarding the equations. Finally, define the space ($4$-groupoid) of adjunctions in a $4$-category $\CC$ to be the full $4$-subgroupoid $\Map^{\adj}(\sk_2(\Adj_{(2,1)}),\CC)\hookrightarrow\Map(\sk_2(\Adj_{(2,1)}),\CC)$ whose objects are those maps $\sk_2(\Adj_{(2,1)})\to\CC$ which determine adjunctions in $\CC$.

Then the question is whether the map $$\Map(\PP,\CC)\to\Map^{\adj}(\sk_2(\Adj_{(2,1)}),\CC)$$ defined by restriction along $\sk_2(\Adj_{(2,1)})\hookrightarrow\PP$ is a trivial fibration. This would mean that, given $(l,r,\un,\co)$ in $\CC$ such that the snake equations are satisfied in $\h_2(\CC)$, the space ($4$-groupoid) of maps $\PP\to\CC$ extending $(l,r,\un,\co)$ would be contractible.
  

This is \textbf{not true} and one can argue that in fact it is \textbf{impossible} to find a computad $\mathcal{Q}$ such that $\Map(\QQ,\CC)\to\Map^{\adj}(\sk_2(\Adj_{(2,1)}),\CC)$ is a trivial fibration. This means it is impossible to find a coherent replacement for the data of an adjunction in $\CC$ as defined above.

So instead we consider the fact that for a $2$-category $\CC$ the data of an adjunction $(l,r,\un,\co)$ is equivalent to the data of a $1$-morphism $l$ in $\CC$ that admits a right adjoint, in the sense that $\Map(\Adj_{(2,1)},\CC)\to\Map^L(\theta^{(1)},\CC)$ is a trivial fibration (see \cite[Proposition 3.3]{adj3} for a proof, although the result was known before). 

With this in mind, given a $4$-category $\CC$, a \textbf{left adjoint $1$-morphism} in $\CC$ consists of
\begin{description}
 \item[Data:] a $1$-morphism $l$ in $\CC$.
 \item[Properties:] there exist $(r,\un,\co)$ in $\CC$ such that $(l,r,\un,\co)$ is an adjunction.
\end{description}

Theorem \ref{intromain} says that this is equivalent to the data of a map $\Adj_{(4,1)}\to\CC$. So in this sense a map $\Adj_{(4,1)}\to \CC$ deserves to called a coherent left adjoint $1$-morphism. However, given the additional results in Section \ref{more_main}, we feel justified in calling it a \textbf{coherent adjunction}, following \cite{riehlverity}.

\subsection{Related work}

The study of coherence for $3$-categorical adjunctions dates back to \cite{verity_thesis}, where the swallowtail relations are first identified in this context. Further developments appear in \cite{gurski}, \cite{bruce} and \cite{piotr}. In \cite{adj3} we proved a coherence theorem for adjunctions in a strict $3$-category, analogous to the one in the present paper.

In \cite{data_struct_quasistrict}, there is a string diagram proof of the fact that, given an adjunction in a $4$-category which is already equipped with $4$-morphisms implementing the swallowtail coherence laws, one can modify one of these $4$-morphisms so that two butterfly relations are satisfied. This proof is formalized in the proof assistant Globular (\cite{globular}).

In \cite{riehlverity} the authors construct an $(\infty,2)$-category $\underline{\Adj}$ and prove that the space of functors $\underline{\Adj}\to \Cat_{\infty}$ is equivalent to the space $1$-morphisms in $\Cat_{\infty}$ which admit a right adjoint. They also prove a result analogous to Theorem \ref{thm2}. One can think of $\Adj_{(4,1)}$ as an explicit finite presentation for the homotopy $4$-category of $\underline{\Adj}$. This seems to be the first place in the literature where coherence for adjunctions in a higher category $\CC$ is stated in terms of an equivalence of spaces between the space of morphisms in $\CC$ which admit an adjoint and the space of maps into $\CC$ out of a category consisting of a free adjunction. The \emph{strictly undulating squiggles} used there are also a kind of string diagram calculus.

\subsection{Future work}

We will use the main result in this paper to give a new proof of the coherence Theorem for $3$-dualizable objects in strict symmetric monoidal $3$-categories in the author's PhD Thesis \cite{araujo_thesis}. This new proof uses the machinery of fibrations of $4$-groupoids and the associated long exact sequence on homotopy groups to vastly reduce the amount and the complexity of the necessary string diagram calculations. In fact, all the necessary string diagram calculations are done in \cite{adj3} and the present paper. The final step is essentially just \cite[Chapter 4]{araujo_thesis} and will be the subject of a subsequent paper.

The cobordism hypothesis should allow us to interpret the presentation encoding a coherent $3$-dualizable object as giving a finite presentation of the $3$-dimensional fully extended framed bordism category, as a symmetric monoidal $3$-category. This would however require the coherence result to be extended to weak symmetric monoidal $3$-categories. These results on coherence for adjoints would easily be extended to the case of weak categories, as long as one establishes that our string diagram calculus is valid in that context. This is also the subject of ongoing research. 

Concretely, we are working on a theory of semistrict $n$-categories based on composition of string diagrams (\cite{sesquicat},\cite{sesqui_comp}). The main obstacle will be establishing the equivalence of this theory to standard models of weak $n$-categories. Recent results suggest this might me achievable (\cite[Section 6.2]{griffiths_thesis}).

A subsequent goal of this project is to use the obtained finite presentation of the bordism category to produce an efficient method for computing the values on $3$-framed manifolds of the fully extended framed Topological Field Theory associated to a fusion category in \cite{dtc}.

\section{Background}

We now give some necessary definitions and recall the main results from \cite{fibrations} and \cite{adj3} which we will need in the rest of this paper.

\subsection{Strict $n$-categories}

We think of a strict $n$-category as an algebra over a certain monad $$T_n:\gSet_n\to\gSet_n$$ on the category of $n$-globular sets. This is the monad defined in \cite[Chapter 8]{operads_cats}. Alternatively one can think of a strict $n$-category as a category enriched in strict $(n-1)$-categories with the cartesian product. We use $\CC_k$ to denote the set of $k$-morphisms in the $n$-category $\CC$. We denote the $k$-fold identity morphism on some $m$-morphism $f$ by $\id^{(k)}_f\in\CC_{m+k}$. A \textbf{full $n$-subcategory} of an $n$-category $\CC$ is an $n$-subcategory $\DD\subset\CC$ such that $\DD(x,y)=\CC(x,y)$ as $(n-1)$-categories, for any objects $x,y\in\CC_0$.

\subsection{Equivalences}

In a strict $n$-category, we say that a $k$-morphism $f:x\to y$ is an \textbf{isomorphism} if there exists another $k$-morphism $f:y\to x$ such that $f\circ g=\id_y$ and $g\circ f=\id_x$. We also say that $f$ is \textbf{invertible} and we call $g$ its \textbf{inverse} (one can show that it is unique). However, we are more interested in a weaker version of this, known as \textbf{equivalence}.

\begin{definition}
	
	Let $\CC$ be a strict $n$-category. An $n$-morphism $f:x\to y$ in $\CC$ is an \textbf{equivalence} if it is an isomorphism. When $k<n$, a $k$-morphism $f:x\to y$ in $\CC$ is an \textbf{equivalence} when there is another $k$-morphism $g:y\to x$ and equivalences $f\circ g\to\id_y$ and $g\circ f\to\id_x$ in $\CC$. We say that $x$ is \textbf{equivalent} to $y$, and write $x\simeq y$, if there is an equivalence $x\to y$. When $f:x\to y$ is an equivalence, we also call it \textbf{weakly invertible} and any morphism $g:y\to x$ such that $f\circ g\simeq \id_y$ and $g\circ f\simeq \id_x$ is called a \textbf{weak inverse} to $f$. When $f$ is a $k$-morphism and an equivalence we also call it a \textbf{$k$-equivalence}.
	
\end{definition}

\begin{definition}
	
	An \textbf{$n$-groupoid} is an $n$-category all of whose morphisms are equivalences.
	
\end{definition}

\begin{definition}
 
Given an strict $n$-category $\CC$ and $k\leq n$, we define its \textbf{homotopy $k$-category} $\h_k(\CC)$ as the $k$-category whose $\ell$-morphisms for $\ell\leq k-1$ agree with those in $\CC$ and whose $k$-morphisms are the equivalence classes of $k$-morphisms in $\CC$. The composition is induced by that in $\CC$. 
 
\end{definition}	

Finally, we use the following notion of weak equivalence for functors, which coincides with the one in the folk model structure of \cite{folk}.

\begin{definition}
	
	A functor $F:\CC\to\DD$ between strict $n$-categories is called \textbf{essentially surjective} if for every object $d\in\DD$ there exists an object $\co\in\CC$ and an equivalence $F(\co)\to d$ in $\DD$. A functor $F:\CC\to\DD$ between strict $n$-categories is called a \textbf{weak equivalence} if it is essentially surjective and for all objects $c_1,c_2\in\CC$ the induced functor $\CC(c_1,c_2)\to\DD(F(c_1),F(c_2))$ is a weak equivalence of $(n-1)$-categories.  
	
\end{definition}

\begin{lemma}\label{comp_weak_eq}

If $F:\CC\to\DD$ and $G:\DD\to \EE$ are weak equivalences of strict $n$-categories, then $G\circ F$ is a weak equivalence.
 
\end{lemma}

\begin{definition}
	
	An $n$-groupoid $G$ is called \textbf{weakly contractible} if the map $G\to *$ is a weak equivalence.
	
\end{definition}

Notice that a weakly contractible $n$-groupoid must be nonempty.

\subsection{Homotopy groups}

\begin{definition}

Let $G$ be an $n$-groupoid and $x\in G_0$ an object. Denote by $\Omega_x G$ the $(n-1)$-groupoid $\Hom(x,x)$. This comes equipped with a strictly associative monoidal structure, given by composition in $G$. 

\end{definition}

\begin{definition}

Let $G$ be an $n$-groupoid. We define $\pi_0(G):=G_0/\sim$, where the equivalence relation $\sim$ is equivalence in $G$. Now let $x\in G$ be an object. Define $\pi_0(G,x)$ to be the pointed set $(\pi_0(G),[x])$, were $[x]$ denotes the equivalence class of $x$ in $\pi_0(G)$. Finally, for $1\leq k\leq n$, define $\pi_k(G,x):=\pi_{k-1}(\Omega_xG,\id_x)$ with monoid structure induced by composition.

\end{definition}

Note that, for $k\geq 1$, the monoids $\pi_k(G,x)$ are actually groups and for $k\geq 2$ they are abelian, by an Eckmann-Hilton argument with pasting diagrams. Moreover, given a map of $n$-groupoids $f:A\to B$ and an object $a\in A$ one can also define $\pi_k(f,a):\pi_k(A,a)\to\pi_k(B,f(a))$, making $\pi_k$ into a functor on pointed $n$-groupoids.

\begin{lemma}\cite[Lemma 5.6]{adj3}\label{htpygrps}

Let $f:A\to B$ be a map of $n$-groupoids. Then $f$ is a weak equivalence if and only if the maps $\pi_k(f,a):\pi_k(A,a)\to\pi_k(B,f(a))$ are isomorphisms, for all $a\in A_0$ and for all $k\geq 0$.

\end{lemma}

\begin{definition}

Given an $n$-groupoid $G$ and an integer $k\geq 0$, we say that $G$ is \textbf{$k$-connected} if $\pi_\ell(G,x)=*$ for every $\ell\leq k$ and $x\in G_0$. We usually write \textbf{connected} when we mean \textbf{$0$-connected}.

\end{definition}

\begin{corollary}\label{nconnected}

An $n$-groupoid $G$ is $n$-connected if and only if it is weakly contractible.
 
\end{corollary}

\subsection{Fibrations}

\begin{definition}
	
	A map of $n$-groupoids $p:E \to B$ is called a \textbf{fibration} if, given any $k$-morphism $f:x\to y$ in $B$ and a lift $\tilde{x}$ of its source along $p$, there exists a lift $\tilde{f}:\tilde{x}\to\tilde{y}$ of $f$ along $p$. A map that is both a fibration and a weak equivalence is called a \textbf{trivial fibration}.
	
\end{definition}

In other words, a fibration of $n$-groupoids is a map that has the right lifting property with respect to source inclusions $s:\theta^{(k-1)}\hookrightarrow\theta^{(k)}$ for all $1\leq k\leq n$.

\begin{remark}\label{folk_fib}
	
	Given $n$-groupoids $E$ and $B$, it is natural to ask whether a map $p:E\to B$ is a fibration in the sense of this paper if and only if it is is a fibration in the folk model structure on strict $n$-categories defined in \cite{folk}. This seems plausible, see \cite{fibrations} for some discussion.
	
\end{remark}

\begin{lemma}\label{comp_fib}

If $f:X\to Y$ and $g:Y\to Z$ are fibrations of strict $n$-groupoids, then $g\circ f$ is a fibration.
 
\end{lemma}

One can define the pullback of a diagram of $n$-categories by taking the pullback in $n$-globular sets and equipping it with a canonical $T_n$-algebra structure. In the case of $n$-groupoids, we have the following result.

\begin{proposition}\cite[Proposition 4.11]{adj3}\label{pbngrpd}
	
	Given a diagram of $n$-groupoids $$\xymatrix{ & X\ar[d]^{F} \\ Y\ar[r]_{G} & Z}$$ where $F$ is a fibration, the pullback $X\times_{Z}Y$ is an $n$-groupoid.
	
\end{proposition}

\begin{definition}
 
Given a fibration of $n$-groupoids $p:E\to B$ and an object $b\in B_0$, its \textbf{fibre} is the $n$-groupoid $p^{-1}(b)$ defined by the pullback \[\begin{tikzcd}
	{p^{-1}(b)} & E \\
	{*} & B
	\arrow[from=1-1, to=1-2]
	\arrow[from=1-1, to=2-1]
	\arrow["\lrcorner"{anchor=center, pos=0.125}, draw=none, from=1-1, to=2-2]
	\arrow["p", from=1-2, to=2-2]
	\arrow["b"', from=2-1, to=2-2]
\end{tikzcd}.\] 
 
\end{definition}

In \cite{adj3} we constructed a long exact sequence in homotopy groups associated to a fibration of $n$-groupoids, which allowed us to prove the following result.

\begin{proposition}\cite[Corollary 5.12]{adj3}\label{fibre}
	
	A fibration of $n$-groupoids $p:E\to B$ is a trivial fibration if and only if for every object $b\in B$ the fibre $p^{-1}(b)$ is weakly contractible.
	
\end{proposition}	

Notice that this implies that such a map must be surjective on objects, as its fibres are contractible, which implies they are nonempty.

We will also need the following standard fact, whose straightforward proof can be found in \cite{adj3}.

\begin{lemma}\label{triv_fib}
	
	Consider a pullback
	
	$$\xymatrix{X\times_B Y\ar[r]\ar[d] & Y \ar[d]^g \\ X\ar[r]_-{f} & B }$$ where $f$ is a trivial fibration. Then $\pi_2:X\times_B Y \to Y$ is a trivial fibration.
	
\end{lemma}

\subsection{Adjunctions}

\begin{definition}
	
	An \textbf{adjunction} in a strict $2$-category $\CC$ is a pair of $1$-morphisms $l:X\to Y$ and $r:Y\to X$ together with $2$-morphisms $\un:\id_X\to r\circ l$ and $\co:l\circ r\to\id_Y$ called the unit and the counit, which satisfy two standard relations, called zigzag, snake or triangle identities. 
	
\end{definition}

The following definitions of adjunctions in $n$-categories are adapted from the ones given in \cite{lurie} for the case of $(\infty,n)$-categories.

\begin{definition}
	
	An \textbf{adjunction} between $1$-morphisms in a strict $n$-category $\CC$ is a pair of $1$-morphisms $l:X\to Y$ and $r:Y\to X$ together with $2$-morphisms $\un:\id_X\to r\circ l$ and $\co:l\circ r\to\id_Y$ called the unit and the counit, which determine an adjunction in the homotopy $2$-category $\h_2(\CC)$.
	
\end{definition}

This means that an adjunction between $1$-morphisms in a $4$-category consists of a pair of $1$-morphisms $l:X\to Y$ and $r:Y\to X$ together with unit and counit $2$-morphisms satisfying the usual snake relations (also known as triangle identities) up to $3$-equivalence.

\begin{definition}
	
	An adjunction between $k$-morphisms in a strict $n$-category $\CC$ is an adjunction between $1$-morphisms in an appropriate $(n-k+1)$-category of morphisms in $\CC$.
	
\end{definition}

The following Lemma relating equivalences and adjunctions is well known (see \cite[Lemma 2.9]{adj3} for a proof).

\begin{lemma}\label{adjeq}
	
	Let $\CC$ be a strict $n$-category, $f:x\to y$ a $k$-equivalence in $\CC$, $g:y\to x$ a weak inverse and $\un:\id_x\to g\circ f$ a $(k+1)$-equivalence. Then there exists a $(k+1)$-equivalence $\co:f\circ g\to \id_y$ such that $(f,g,\un,\co)$ is an adjunction in $\CC$.
	
\end{lemma}

\subsection{Presentations}

An \textbf{$n$-categorical presentation} is simply a collection of $k$-cells for every $k\leq n+1$, whose sources and targets are composites of lower dimensional cells. We interpret the $(n+1)$-cells as relations. Given an $n$-categorical presentation $\PP$ we denote by $F(\PP)$ the \textbf{$n$-category generated by $\PP$}. Its $k$-morphisms are arbitrary composites of the $k$-cells in $\PP$. Two $n$-morphisms are declared equal when they are related by an $(n+1)$-cell. We sometimes write $\PP\to\CC$ to refer to a functor $F(\PP)\to\CC$.

This can be made precise by using the theory of \textbf{computads}. See \cite{csp_thesis} for a detailed treatment of computads.

We denote by $\theta^{(k)}$ the computad generated by a single $k$-cell, so that functors $\theta^{(k)}\to\CC$ are in canonical bijection with the set of $k$-morphisms in $\CC$. 

\subsection{String diagrams}

We refer the reader to \cite{fibrations} for a detailed description of the string diagram calculus for $4$-categories which we will use in this paper. However, we hope the reader who is familiar with the general idea of string diagrams for $2$-categories (or monoidal categories) may be able to read on without reading \cite{fibrations} in detail.

\subsection{Functor categories}

In \cite{fibrations} we gave an explicit description of the $4$-category $\Fun(\CC,\DD)$ in terms of string diagrams, when $\CC$ and $\DD$ are $4$-categories. We include it below for convenience. We denote by $\Map(\CC,\DD)$ the underlying $4$-groupoid in $\Fun(\CC,\DD)$. Given a presentation $\PP$ we write $\Fun(\PP,\DD)$ instead of $\Fun(F(\PP),\DD)$ and similarly for $\Map$. We denote by $\Map^L(\theta^{(1)},\CC)$ the full $4$-subgroupoid of $\Map(\theta^{(1)},\CC)$ whose objects are the $1$-morphisms in $\CC$ which admit a right adjoint in $\h_2(\CC)$.

\subsubsection{Natural transformations}

Given functors $F,G:\CC\to\DD$, a \textbf{natural transformation}, or $1$-transfor, $\alpha:F\to G$ consists of the following data. We use {\red} and {\blue} to denote the images of objects and morphisms under $F$
and $G$, respectively.
\begin{enumerate}
	\item[0.] For each object $Y\in \CC$ a $1$-morphism $\alpha_Y:F(Y)\to G(Y)$. \begin{center}\begin{tabular}{lcr}$Y=$ \includegraphics[scale=1.5,align=c]{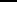} & $\mapsto$ & $\alpha_Y=$  \includegraphics[scale=1.5,align=c]{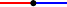}\end{tabular}\end{center}
	
	\item For each $1$-morphism $g:X\to Y$ in $\CC$ an invertible $2$-morphism $\alpha_g$ in $\DD$. \begin{center}\begin{tabular}{lcr} $g=$ \includegraphics[scale=1.5,align=c]{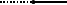} & $\mapsto$ & $\alpha_g=$ \includegraphics[scale=1.5,align=c]{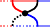} $:$ \includegraphics[scale=1.5,align=c]{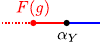} $\to$ \includegraphics[scale=1.5,align=c]{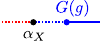}\end{tabular}\end{center}
	
	\item For each $2$-morphism $\zeta:f\to g$ in $\CC$ an invertible $3$-morphism $\alpha_{\zeta}$ in $\DD$. \begin{center}\begin{tabular}{lcr} $\zeta=$ \includegraphics[scale=1.5,align=c]{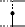} & $\mapsto$ & $\alpha_{\zeta}=$ \includegraphics[scale=1.5,align=c]{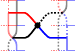} $:$ \includegraphics[scale=1.5,align=c]{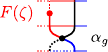} $\to$ \includegraphics[scale=1.5,align=c]{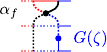}\end{tabular} \end{center}
	
	\item For each $3$-morphism $t:\eta\to\zeta$ in $\CC$ an invertible $4$-morphism $\alpha_{t}$ in $\DD$. \begin{center}\begin{tabular}{lcr}$t=$ \includegraphics[scale=1.5,align=c]{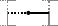} & $\mapsto$ & $\alpha_{t}=$ \includegraphics[scale=1.5,align=c]{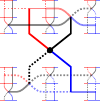} $:$ \includegraphics[scale=1.5,align=c]{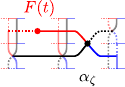} $\to$ \includegraphics[scale=1.5,align=c]{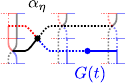}\end{tabular}\end{center}
	
	\item For each $4$-morphism $\W:s\to t$ in $\CC$ a relation $\alpha_{\W}$ in $\DD$.  \begin{center}\begin{tabular}{lcr} $\W=$ \includegraphics[scale=1.5,align=c]{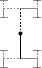} & $\mapsto$ & $\alpha_{\W}:$\includegraphics[scale=1.5,align=c]{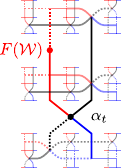} $=$ \includegraphics[scale=1.5,align=c]{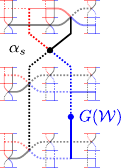}\end{tabular}\end{center}
	
\end{enumerate}

This data is subject to relations equating the values of $\alpha$ on composite morphisms with the corresponding composites of values of $\alpha$ given by stacking diagrams.

\subsubsection{Modifications}

Given natural transformations $\alpha,\beta:F\to G$, a \textbf{modification}, or $2$-transfor, $m:\alpha\to\beta$ consists of the following data. We use {\green} for $\alpha$ and {\purple} for $\beta$. 

\begin{enumerate}
	\item[0.] For each object $Y\in\CC$ a $2$-morphism $m_Y:\alpha_Y\to \beta_Y$ in $\DD$. \begin{center}\begin{tabular}{lcr} $Y=$ \includegraphics[scale=1.5,align=c]{functorcat/1morph/y.pdf} & $\mapsto$ & $m_Y=$ \includegraphics[scale=1.5,align=c]{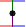} $:$ \includegraphics[scale=1.5,align=c]{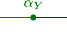} $\to$ \includegraphics[scale=1.5,align=c]{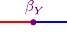}\end{tabular}\end{center}
	
	\item For each $1$-morphism $g:X\to Y$ in $\CC$ an invertible $3$-morphism $m_g$ in $\DD$. \begin{center}\begin{tabular}{lcr}$g=$ \includegraphics[scale=1.5,align=c]{morphisms/f.pdf} & $\mapsto$ & $m_g=$ \includegraphics[scale=1.5,align=c]{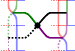} $:$ \includegraphics[scale=1.5,align=c]{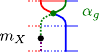} $\to$ \includegraphics[scale=1.5,align=c]{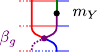}\end{tabular}\end{center}
	
	\item For each $2$-morphism $\zeta=\includegraphics[scale=1.5,align=c]{morphisms/eta.pdf}:f\to g$ in $\CC$ an invertible $4$-morphism $m_{\zeta}$ in $\DD$.  $$m_{\zeta}= \includegraphics[scale=1.5,align=c]{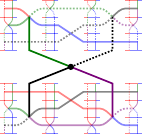}: \includegraphics[scale=1.5,align=c]{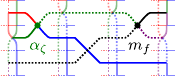}\to \includegraphics[scale=1.5,align=c]{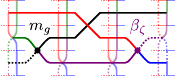}$$
	
	\item For each $3$-morphism $t:\eta\to\zeta$ in $\CC$ a relation $m_t$ in $\DD$. \begin{center}\begin{tabular}{lcr}$t=$ \includegraphics[scale=1.5,align=c]{morphisms/s.pdf} & $\mapsto$ & $m_{t}:$ \includegraphics[scale=1.5,align=c]{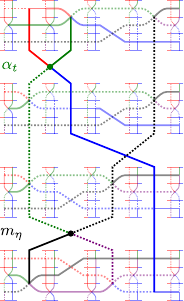} $=$ \includegraphics[scale=1.5,align=c]{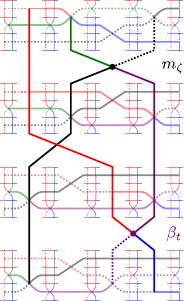}\end{tabular}\end{center}
	
\end{enumerate}

This data is subject to relations equating the values of $m$ on composite morphisms with the corresponding composites of values of $m$ given by stacking diagrams.

\subsubsection{Perturbations}

Given modifications $l,m:\alpha\to\beta$, a \textbf{perturbation}, or $3$-transfor, $\A:l\to m$ consists of the following data. We use {\orange} for $l$ and {\lightblue} for $m$. 

\begin{enumerate}
	\item[0.] For each object $Y\in\CC$ a $3$-morphism $\A_Y:l_Y\to m_Y$ in $\DD$. \begin{center}\begin{tabular}{lcr} $Y=$ \includegraphics[scale=1.5,align=c]{functorcat/1morph/y.pdf} & $\mapsto$ & $\A_Y=$\includegraphics[scale=1.5,align=c]{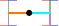} $:$ \includegraphics[scale=1.5,align=c]{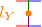} $\to$ \includegraphics[scale=1.5,align=c]{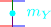}\end{tabular}\end{center}
	
	\item For each $1$-morphism $g=\includegraphics[scale=1.2,align=c]{morphisms/f.pdf}:X\to Y$ in $\CC$ an invertible $4$-morphism $\A_g$ in $\DD$. \begin{center} $\A_g=$ \includegraphics[scale=1.5,align=c]{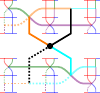} $:$ \includegraphics[scale=1.5,align=c]{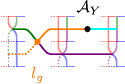} $\to$ \includegraphics[scale=1.5,align=c]{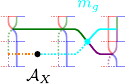}\end{center}
	
	\item For each $2$-morphism $\zeta:f\to g$ in $\CC$ a relation $\A_{\zeta}$ in $\DD$. \begin{center}\begin{tabular}{lcr} $\zeta=$ \includegraphics[scale=1.5,align=c]{morphisms/eta.pdf} & $\mapsto$ & $\A_{\zeta}:$ \includegraphics[scale=1.5,align=c]{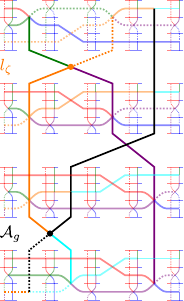} $=$ \includegraphics[scale=1.5,align=c]{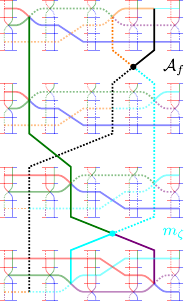}\end{tabular}\end{center}
	
\end{enumerate}

This data is subject to relations equating the values of $\A$ on composite morphisms with the corresponding composites of values of $\A$ given by stacking diagrams. 

\subsubsection{$4$-transfors}

Given perturbations $\A,\B:l\to m$, a \textbf{$4$-transfor} $\Lambda:\A\to\B$ consists of the following data. We use {\olive} for $\A$ and {\teal} for $\B$.  

\begin{enumerate}
	\item[0.] For each object $Y\in\CC$ a $4$-morphism $\Lambda_Y:\A_Y\to \B_Y$ in $\DD$.  \begin{center}\begin{tabular}{lcr} $Y=$ \includegraphics[scale=1.5,align=c]{functorcat/1morph/y.pdf} & $\mapsto$ & $\Lambda_Y=$ \includegraphics[scale=1.5,align=c]{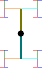} $:$ \includegraphics[scale=1.5,align=c]{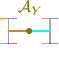} $\to$ \includegraphics[scale=1.5,align=c]{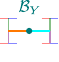}\end{tabular}\end{center}
	
	\item For each $1$-morphism $g=:X\to Y$ in $\CC$ a relation $\Lambda_g$ in $\DD$. \begin{center}\begin{tabular}{lcr}$g=$  \includegraphics[scale=1.5,align=c]{morphisms/f.pdf} & $\mapsto$ & $\Lambda_g:$ \includegraphics[scale=1.5,align=c]{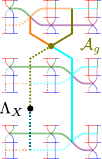} $=$ \includegraphics[scale=1.5,align=c]{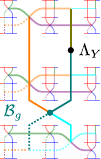}\end{tabular}\end{center}
	
\end{enumerate}

\subsection{Fibrations of mapping $4$-groupoids}

\begin{theorem}\cite[Theorem 1.2]{fibrations}\label{fibration}
	
	Let $\CC$ be a strict $4$-category, $\PP$ a presentation and $\QQ$ another presentation, obtained by adding a finite number of cells to $\PP$. Then the restriction map $$\Map(\QQ,\CC)\to\Map(\PP,\CC)$$ is a fibration of $4$-groupoids.
	
\end{theorem} 	

In this paper we will need the following strenghtening of Theorem \ref{fibration}.

\begin{corollary}\label{corfib}

Let $\CC$ be a strict $4$-category, $\PP$ a presentation and $\QQ$ another presentation, obtained by adding a finite number of cells to $\PP$. Then the restriction map $$\Map(\QQ,\CC)\to\Map(\QQ,\h_k(\CC))\times_{\Map(\PP,\h_k(\CC))}\Map(\PP,\CC)$$ is a fibration of $4$-groupoids. 
 
\end{corollary}

First notice that the target is indeed a $4$-groupoid, because $\Map(\QQ,\h_k(\CC))\to\Map(\PP,\h_k(\CC))$ is a fibration, by Theorem \ref{fibration}. 

\begin{proposition}\label{basefib}

Let $\CC$ a strict $4$-category, $\PP$ a presentation, to which we add an $m$-cell $g$. Then the restriction map $$\Map(\PP\cup\{g\},\CC)\to\Map(\PP\cup\{g\},\h_k(\CC))\times_{\Map(\PP,\h_k(\CC))}\Map(\PP,\CC)$$ is a fibration of $4$-groupoids. 
 
\end{proposition}

\begin{proof}

Suppose we have an $\ell$-morphism $\alpha:F\to G$ in $\Map(\PP\cup\{g\},\h_k(\CC))\times_{\Map(\PP,\h_k(\CC))}\Map(\PP,\CC)$ and a lift $\tilde{F}$ of $F$ to $\Map(\PP\cup\{g\},\CC)$. We want to find a lift $\tilde{\alpha}:\tilde{F}\to\tilde{G}$ in $\Map(\PP\cup\{g\},\CC)$. So we have $\alpha_g\in\h_k(\CC)_{m+\ell}$, $G_g\in\h_k(\CC)_{m+\ell-1}$ and we want to lift both to $\CC$.

If $m+\ell\leq k+1$ then one can easily find the desired lifts, simply by the definition of $\h_k(\CC)$. The reader who is interested in checking the details can split further into cases $m+\ell<k$, $m+\ell=k$ and $m+\ell=k+1$.

If $m+\ell\geq k+2$ then the given data $\alpha_g\in\h_k(\CC)_{m+\ell}$ and $G_g\in\h_k(\CC)_{m+\ell-1}$ is actually no data at all, so this becomes a lifting problem for $\Map(\PP\cup\{g\},\CC)\to\Map(\PP,\CC)$, which can be solved by Theorem \ref{fibration}.\end{proof}

\begin{proof}(of Corollary  \ref{corfib})
 
We use induction on the number of cells added to $\PP$ to obtain $\QQ$. The base case is Proposition \ref{basefib}. For the induction step, suppose $g$ is the last cell added in getting from $\PP$ to $\QQ$. Then the restriction map $$\Map(\QQ,\CC)\to\Map(\QQ,\h_k(\CC))\times_{\Map(\PP,\h_k(\CC))}\Map(\PP,\CC)$$ is equal to the composite \[\begin{tikzcd}
	{\Map(\QQ,\CC)} \\
	{\Map(\QQ,\h_k(\CC))\times_{\Map(\QQ\setminus\{g\},\h_k(\CC))}\Map(\QQ\setminus\{g\},\CC)} \\
	{\Map(\QQ,\h_k(\CC))\times_{\Map(\QQ\setminus\{g\},\h_k(\CC))}(\Map(\QQ\setminus\{g\},\h_k(\CC))\times_{\Map(\PP,\h_k(\CC))}\Map(\PP,\CC))} \\
	{\Map(\QQ,\h_k(\CC))\times_{\Map(\PP,h_k(\CC))}\Map(\PP,\CC)}
	\arrow["fib", from=1-1, to=2-1]
	\arrow["{\id\times_{\id} fib}", from=2-1, to=3-1]
	\arrow["\simeq", from=3-1, to=4-1]
\end{tikzcd}\] where the first map is a fibration by Proposition \ref{basefib}, the second map is a fibration by the induction hypothesis and the final map is a fibration because it is an isomorphism. Finally, the composite is a fibration by Lemma \ref{comp_fib}.\end{proof}

\subsection{Coherence for adjunctions in a 3-category}

We recall here the coherence result for adjunctions in a $3$-category from \cite{adj3}, which we will need to use. We defined the $3$-categorical presentation $\Adj_{(3,1)}$ as follows.

\begin{definition}\cite[Definition 6.2]{adj3}
	
	The presentation $\Adj_{(3,1)}$ consists of 
	
	\begin{enumerate}[start=0,label={(\arabic*)}]
		\item objects $X=$ \includegraphics[scale=2,align=c]{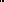} and $Y=$ \includegraphics[scale=2,align=c]{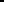} 
		\item $1$-cells $l$= \includegraphics[scale=1,align=c]{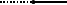} $:\xymatrix@1{X\ar@1@<1ex>[r] & \ar@1@<1ex>[l]Y}:$ \includegraphics[scale=1,align=c]{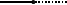} $=r$
		\item $2$-cells

		\begin{center}\begin{tabular}{llllll}

				$\un=$ & \includegraphics[scale=1,align=c]{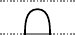} & $:$ & \includegraphics[scale=1,align=c]{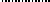} & $\xymatrix@1{\ar@2[r] & }$ & \includegraphics[scale=1,align=c]{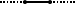} 
				
				\\
				
				\\
				
				$\co=$ & \includegraphics[scale=1,align=c]{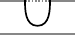} & $:$ & \includegraphics[scale=1,align=c]{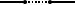} & $\xymatrix@1{\ar@2[r] & }$ & \includegraphics[scale=1,align=c]{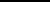} 
				
			\end{tabular}\end{center}

			\item $3$-cells
			
			\begin{center}\begin{tabular}{lllllllll}

					$C_l=$ & \includegraphics[scale=1,align=c]{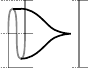} & $:$ & \includegraphics[scale=1,align=c]{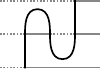} & $\xymatrix@1{\ar@3@<1ex>[r] & \ar@3@<1ex>[l]}$ & \includegraphics[scale=1,align=c]{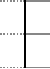} & $:$ & \includegraphics[scale=1,align=c]{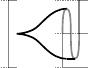} & $=C_l^{-1}$
					
					\\
					
					\\
					
					$C_r=$ & \includegraphics[scale=1,align=c]{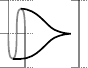} & $:$ & \includegraphics[scale=1,align=c]{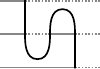} & $\xymatrix@1{\ar@3@<1ex>[r] & \ar@3@<1ex>[l]}$ & \includegraphics[scale=1,align=c]{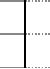} & $:$ & \includegraphics[scale=1,align=c]{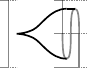} & $=C_r^{-1}$
					
				\end{tabular}\end{center}

				\item relations
				
				\begin{center}\begin{tabular}{lllll}

						\includegraphics[scale=1,align=c]{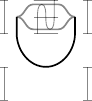} & $:$ & \includegraphics[scale=1,align=c]{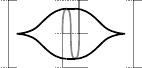} & $=$ & $\Id^{(2)}_l$ 
						
						\\
						
						\\

						\includegraphics[scale=1,align=c]{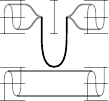} & $:$ & \includegraphics[scale=1,align=c]{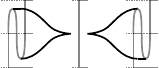} & $=$ & \includegraphics[scale=1,align=c]{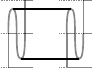} 
						
						\\
						
						\\
						
						\includegraphics[scale=1,align=c]{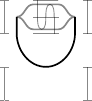} & $:$ & \includegraphics[scale=1,align=c]{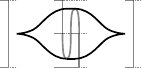} & $=$ & $\Id^{(2)}_r$ 
						
						\\
						
						\\

						\includegraphics[scale=1,align=c]{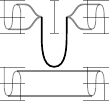} & $:$ & \includegraphics[scale=1,align=c]{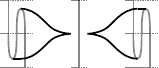} & $=$ & \includegraphics[scale=1,align=c]{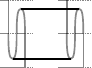}

						\\

						\\

						\includegraphics[scale=1,align=c]{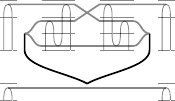} & $:$ & \includegraphics[scale=1,align=c]{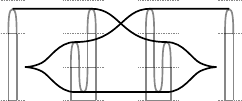} & $=$ & \includegraphics[scale=1,align=c]{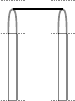}
						
					\end{tabular}\end{center}

				\end{enumerate}
				
			\end{definition}
						
We proved the corresponding coherence result.
			
\begin{theorem}\cite[Theorem 1.1]{adj3}\label{adj3}
	
	Given a strict $3$-category $\CC$, the restriction map $$E_l:\Map(\Adj_{(3,1)},\CC)\to\Map^{L}(\theta^{(1)}, \CC)$$ is a weak equivalence of strict $3$-groupoids.
	
\end{theorem}		

Finally, we showed that certain additional relations are satisfied in $F(\Adj_{(3,1)})$. Denote by $(\SW)$ the last relation in $\Adj_{(3,1)}$, which is called a \textbf{swallowtail} relation. Usually the definition of a coherent adjunction in a $3$-category includes another swallowtail relation, which we call here $(\SW_2)$, but in \cite{piotr} the author proves that one follows from the other. In \cite{adj3} we gave a new string diagram proof of this fact. There are also relations $(\overline{\SW})$ and $(\overline{\SW_2})$ whose sources are inverse to those of $\SW$ and $\SW_2$. We therefore defined a larger presentation $\Adj^+_{(3,1)}$ and showed that $F(\Adj_{(3,1)}^{+})=F(\Adj_{(3,1)})$.

\begin{definition}\cite[Definition 6.9]{adj3}
	
	We define $\Adj^{+}_{(3,1)}$ to be the presentation obtained from $\Adj_{(3,1)}$ by adding the relations
	
	\begin{center}
		\begin{center}\begin{tabular}{llllll}

				$(\overline{\SW})$ & \includegraphics[scale=1,align=c]{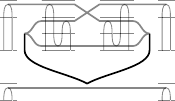} & $:$ & \includegraphics[scale=1,align=c]{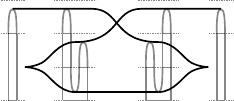} & $=$ & \includegraphics[scale=1,align=c]{adj/pres/id_u.pdf} 
				
				\\
				
				\\
				
				$(\SW_2)$ & \includegraphics[scale=1,align=c]{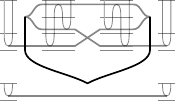} & $:$ & \includegraphics[scale=1,align=c]{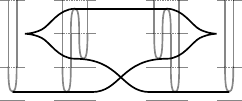} & $=$ & \includegraphics[scale=1,align=c]{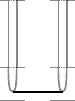} 
				
				\\
				
				\\
				
				$(\overline{\SW_2})$ & \includegraphics[scale=1,align=c]{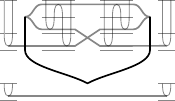} & $:$ & \includegraphics[scale=1,align=c]{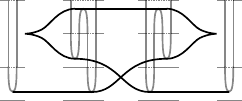} & $=$ & \includegraphics[scale=1,align=c]{adj/pres/id_c.pdf}. 
				
			\end{tabular}\end{center}
			
		\end{center}
		
	\end{definition}
	
	\begin{proposition}\cite[Proposition 6.10]{adj3}\label{sw2}
		
		We have $F(\Adj_{(3,1)}^{+})=F(\Adj_{(3,1)})$. 
		
	\end{proposition}

\section{The presentation}
 
We start by defining the presentation $\Adj_{(4,1)}$ consisting of coherence data for an adjunction between $1$-morphisms in a $4$-category. This presentation is obtained from $\Adj_{(3,1)}$ by replacing each relation between $3$-morphisms by a pair  of $4$-morphisms and then adding relations witnessing the fact that these $4$-morphisms are inverse to each other. Moreover, we add two snake relations expressing the fact that each of the two weakly inverse pairs of cusp $3$-cells in $\Adj_{(4,1)}$ is an adjoint equivalence (note that two adjunctions would normally require four snake relations, but in the case of adjoint equivalences one snake relation implies the other, see the nLab page "adjoint equivalence" for a string diagram proof of this fact).
 
\begin{definition}
 
The presentation $\Adj_{(4,1)}$ consists of 

\begin{enumerate}[start=0,label={(\arabic*)}]
 \item $0$-cells $X=$ \includegraphics[align=c,scale=2]{adj/pres/x.pdf} and $Y=$ \includegraphics[align=c,scale=2]{adj/pres/y.pdf} 
 \item $1$-cells $$l= \includegraphics[align=c,scale=1]{adj/pres/l.pdf} :\xymatrix@1{X\ar@1@<1ex>[r] & \ar@1@<1ex>[l]Y}: \includegraphics[align=c,scale=1]{adj/pres/r.pdf} =r$$
 \item $2$-cells 
 
 \begin{center}
  \begin{tabular}{cccccc}

 $\un=$ & \includegraphics[align=c,scale=1]{adj/pres/u.pdf} & $:$ & \includegraphics[align=c,scale=1]{adj/pres/u_s.pdf} & $\xymatrix@1{\ar@2[r] & }$ & \includegraphics[align=c,scale=1]{adj/pres/u_t.pdf} \\ \\
 
  $\co=$ & \includegraphics[align=c,scale=1]{adj/pres/c.pdf} & $:$ & \includegraphics[align=c,scale=1]{adj/pres/c_s.pdf} & $\xymatrix@1{\ar@2[r] & }$ & \includegraphics[align=c,scale=1]{adj/pres/c_t.pdf} 
 
   \end{tabular}

 \end{center}
 
 \item $3$-cells
 
  \begin{center}
  \begin{longtable}{ccccccccc}

 $C_l=$ & \includegraphics[align=c,scale=1]{adj/pres/cusp_l.pdf} & $:$ & \includegraphics[align=c,scale=1]{adj/pres/snake_l.pdf} & $\xymatrix@1{\ar@3@<1ex>[r] & \ar@3@<1ex>[l]}$ & \includegraphics[align=c,scale=1]{adj/pres/id_l.pdf} & $:$ & \includegraphics[align=c,scale=1]{adj/pres/cusp_l_inv.pdf} & $=C_l^{-1}$
 
 \\
 
 \\
 
 $C_r=$ & \includegraphics[align=c,scale=1]{adj/pres/cusp_r.pdf} & $:$ & \includegraphics[align=c,scale=1]{adj/pres/snake_r.pdf} & $\xymatrix@1{\ar@3@<1ex>[r] & \ar@3@<1ex>[l]}$ & \includegraphics[align=c,scale=1]{adj/pres/id_r.pdf} & $:$ & \includegraphics[align=c,scale=1]{adj/pres/cusp_r_inv.pdf} & $=C_r^{-1}$
 
   \end{longtable}

 \end{center}
 
 \item $4$-cells
 
 \begin{center}
  \begin{longtable}{ccccccccc}

 $c_{C_l}=$ & \includegraphics[align=c,scale=1]{adj/pres/cusp_l_cancel1.pdf}  & $:$ & \includegraphics[align=c,scale=0.7]{adj/pres/cusp_l_cancel1_s.pdf} & $\xymatrix@1{\ar@1@<1ex>[r] & \ar@1@<1ex>[l]}$ & $\Id^{(2)}_l$ & $:$ &  \includegraphics[align=c,scale=1]{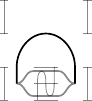} & $=c_{C_l}^{-1}$ 
 
 \\
 
 \\

 $u_{C_l}^{-1}=$ & \includegraphics[align=c,scale=1]{adj/pres/cusp_l_cancel2.pdf}  & $:$ & \includegraphics[align=c,scale=0.7]{adj/pres/cusp_l_cancel2_s.pdf} & $\xymatrix@1{\ar@1@<1ex>[r] & \ar@1@<1ex>[l]}$ & \includegraphics[align=c,scale=0.7]{adj/pres/id_snake_l.pdf} & $:$ &  \includegraphics[align=c,scale=1]{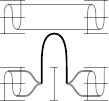} & $=u_{C_l}$ 
 
  \\
 
 \\
 
  $c_{C_r}=$ & \includegraphics[align=c,scale=1]{adj/pres/cusp_r_cancel1.pdf}  & $:$ & \includegraphics[align=c,scale=0.7]{adj/pres/cusp_r_cancel1_s.pdf} & $\xymatrix@1{\ar@1@<1ex>[r] & \ar@1@<1ex>[l]}$ & $\Id^{(2)}_r$ & $:$ &  \includegraphics[align=c,scale=1]{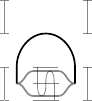} & $=c_{C_r}^{-1}$
 
 \\
 
 \\

 $u_{C_r}^{-1}=$ & \includegraphics[align=c,scale=1]{adj/pres/cusp_r_cancel2.pdf}  & $:$ & \includegraphics[align=c,scale=0.7]{adj/pres/cusp_r_cancel2_s.pdf} & $\xymatrix@1{\ar@1@<1ex>[r] & \ar@1@<1ex>[l]}$ & \includegraphics[align=c,scale=0.7]{adj/pres/id_snake_r.pdf} & $:$ &  \includegraphics[align=c,scale=1]{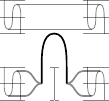}& $=u_{C_r}$

  \\

 \\

  $\SW=$ &\includegraphics[align=c,scale=1]{adj/pres/swallowtail.pdf}  & $:$ & \includegraphics[align=c,scale=0.7]{adj/pres/swallowtail_s.pdf} & $\xymatrix@1{\ar@1@<1ex>[r] & \ar@1@<1ex>[l]}$ & \includegraphics[align=c,scale=0.7]{adj/pres/id_u.pdf} &  $:$ &  \includegraphics[align=c,scale=1]{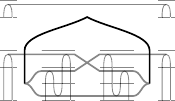} & $=SW^{-1}$
 
   \end{longtable}

 \end{center}
 
 \item relations
 
 \begin{center}
  \begin{longtable}{ccccccc}
   
\includegraphics[align=c,scale=1]{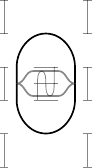}  & $=$ & $\Id^{(3)}_l$ & ; & \includegraphics[align=c,scale=1]{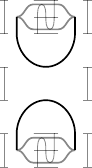}  & $=$ & \includegraphics[align=c,scale=1]{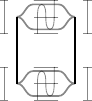}

\\

\\

\includegraphics[align=c,scale=1]{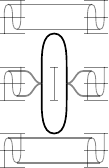}  & $=$ & \includegraphics[align=c,scale=1]{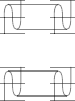} & ; & \includegraphics[align=c,scale=1]{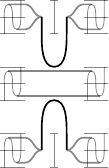}  & $=$ & \includegraphics[align=c,scale=1]{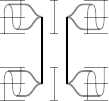}

\\

\\

\includegraphics[align=c,scale=1]{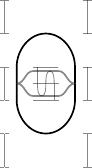}  & $=$ & $\Id^{(3)}_r$ & ; & \includegraphics[align=c,scale=1]{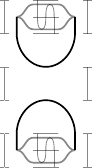}  & $=$ & \includegraphics[align=c,scale=1]{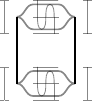}

\\

\\

\includegraphics[align=c,scale=1]{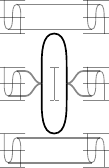}  & $=$ & \includegraphics[align=c,scale=1]{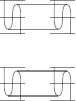} & ; & \includegraphics[align=c,scale=1]{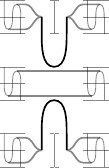}  & $=$ & \includegraphics[align=c,scale=1]{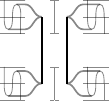}

\\
\\

\includegraphics[align=c,scale=1]{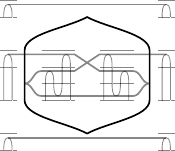}  & $=$ & $\Id^{(2)}_u$ & ; & \includegraphics[align=c,scale=1]{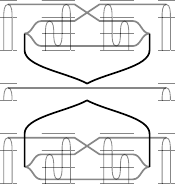}  & $=$ & \includegraphics[align=c,scale=1]{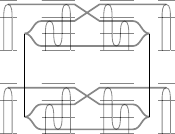}

\\
\\

\includegraphics[align=c,scale=1]{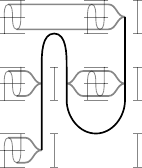}  & $=$ & $\Id_{C_l}$ & ; & \includegraphics[align=c,scale=1]{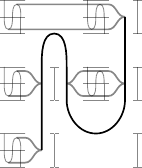}  & $=$ & $\Id_{C_r}$ 

   \end{longtable}

 \end{center}

\end{enumerate}

\end{definition}

An adjunction in a $4$-category is a pair of $1$-morphisms, together with unit and counit $2$-morphisms satisfying the snake relations up to $3$-equivalence. By picking weakly inverse pairs of $3$-morphisms - which are usually called \textbf{cusp $3$-morphisms} - implementing these snake relations and then considering the invertible $4$-morphisms and relations between these, implied in witnessing the fact that these are in fact weakly inverse pairs of $3$-morphisms, one obtains most of the cells in the above presentation. Additionally, one needs the \textbf{swallowtail} $4$-morphism $\SW$ and its inverse, together with two relations witnessing the fact that these are an inverse pair. Finally, we add two snake relations expressing the fact that each weakly inverse pair of cusp $3$-morphisms is an adjoint equivalence.

In \cite{verity_thesis} and \cite{gurski} the definitions of coherent adjunction include two swallowtail relations, but in \cite{piotr} the author shows that the second follows from the first. Here we are one categorical dimension above, so we instead have a swallowtail $3$-morphism (together with an inverse), from which we can construct a second swallowtail $3$-morphism, related to the first by a \textbf{butterfly relation}. For this reason, we only need to include one swallowtail $3$-morphism and no butterfly relations in our presentation. The butterfly relations are included as part of the coherence data for an adjunction in a $4$-category in e.g. \cite{data_struct_quasistrict}. Our result shows that they are actually not required in a minimal presentation of the coherence data. We will show below how to enlarge our presentation to include the butterfly relations, as well as certain relations we call \textbf{swallowtail flip relations} (analogous to the cusp flip relations in \cite{csp_thesis}), obtaining a presentation which generates an isomorphic $4$-category.

\subsection{Enlarging the presentation}

In writing down the presentation $\Adj_{(4,1)}$ we have tried to make it as small as possible. However, it is also sometimes useful to have more generators and relations at our disposal. We now introduce an enlarged presentation $\Adj_{(4,1)}^{max}$ and show that the map $F(\Adj_{(4,1)})\to F(\Adj_{(4,1)}^{max})$ induced by the inclusion of presentations is an isomorpism of $4$-categories, by constructing an explicit inverse. We do this by introducing new generating $4$-cells, together with relations which force them to be equal to certain composites of generators in $\Adj_{(4,1)}$, as well as some extra relations which already follow from the ones in $\Adj_{(4,1)}$. 

\subsubsection{Additional snake relations}

Consider the final two relations in $\Adj_{(4,1)}$. We include them because we want the cusp $3$-morphisms to be not only equivalences but adjoint equivalences. More precisely, we want $(C_l\dashv C_l^{-1},u_{C_l},c_{C_l})$ and $(C_r\dashv C_r^{-1},u_{C_r},c_{C_r})$ to be adjunctions. However, we have only included two snake relations, whereas we would need four such relations to witness these two adjunctions. The reason we left out these other snake relations is that, in the case of adjoint equivalences, they follow from the ones we included. A string diagram proof of this fact can be found in the nLab page "adjoint equivalence". We can also include the four snake relations involving the generators $u_{C_l}^{-1}$, $c_{C_l}^{-1}$, $u_{C_r}^{-1}$, $c_{C_r}^{-1}$, which follow from the ones above because their sources are inverses.

\begin{definition}
	
Let $\Adj_{(4,1)}^+$ be the presentation obtained from $\Adj_{(4,1)}$ by adding the six relations 	\begin{center}\begin{tabular}{ccccccc}\includegraphics[align=c,scale=1]{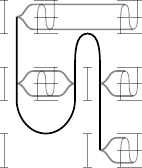}  & $=$ & $\Id_{C_l^{-1}}$ & ; &

\includegraphics[align=c,scale=1]{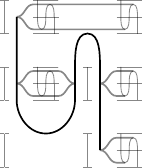}  & $=$ & $\Id_{C_r^{-1}}$

\\

\\

\includegraphics[align=c,scale=1]{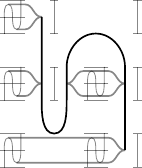} & $=$ & $\Id_{C_l}$ & ; & \includegraphics[align=c,scale=1]{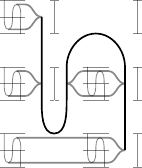} & $=$ & $\Id_{C_r}$  

\\

\\

 \includegraphics[align=c,scale=1]{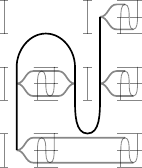} & $=$ & $\Id_{C_l^{-1}}$ & ; & \includegraphics[align=c,scale=1]{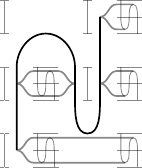} & $=$ & $\Id_{C_r^{-1}}.$
	\end{tabular}\end{center} 

\end{definition}	

\begin{proposition}

We have $F(\Adj_{(4,1)})=F(\Adj_{(4,1)}^+)$. 	
	
\end{proposition}	

\begin{proof}

We have only added relations which already hold in $F(\Adj_{(4,1)})$.\end{proof}

\subsubsection{The second swallowtail $4$-morphism and the butterfly relations}\label{second_swallowtail}

We will now construct inside $F(\Adj_{(4,1)})$ a second swallowtail $4$-morphism $\SW_2$, together with an inverse $\SW_2^{-1}$.

$$\SW_2=\includegraphics[align=c,scale=1]{adj/obj/extra/swallowtail3.pdf} : \includegraphics[align=c,scale=1]{adj/obj/extra/swallowtail3_s.pdf} \xymatrix@1{\ar@1@<1ex>[r] & \ar@1@<1ex>[l]} \includegraphics[align=c,scale=1]{adj/pres/id_c.pdf} : \includegraphics[align=c,scale=1]{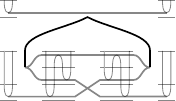}=\SW_2^{-1}.$$ We will then show that these $4$-morphisms are related to $\SW,\SW^{-1}$ by two butterfly relations.

The construction proceeds as follows. Consider the map of $3$-categories $\F(\Adj_{(3,1)})\to \h_3(\F(\Adj_{(4,1)}))$ determined by the obvious map on generators. By Proposition \ref{sw2} the relation $\SW_2$ holds in $\F(\Adj{(3,1)})$, so it will also hold in $\h_3(\F(\Adj_{(4,1)}))$. Therefore we can pick an isomorphism 

\begin{center}\includegraphics[align=c,scale=1]{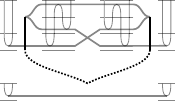} $:$ \includegraphics[align=c,scale=1]{adj/obj/extra/swallowtail3_s.pdf} $\xymatrix@1{\ar@1[r]&}$ \includegraphics[align=c,scale=1]{adj/pres/id_c.pdf}\end{center} in $\F(\Adj_{(4,1)})$. This is obviously not unique, but an explicit construction can be found in \cite{adj3}. Now we just need to modify this so that it satisfies the butterfly relations. So define $\SW_2$ by

\begin{center}
	
	$\SW_2=$\includegraphics[align=c,scale=1]{adj/obj/extra/swallowtail3.pdf} $:=$ \includegraphics[align=c,scale=1]{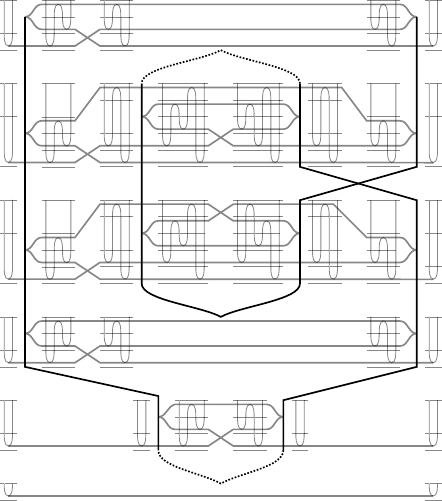}.  
	
\end{center}

It is clear that $\SW_2$ is invertible, being a composite of invertible morphisms. 

\begin{lemma}\label{butterfly1}
	
The $4$-morphisms $\SW^{-1}$ and $\SW_2$ (as defined above) satisfy the butterfly relation 	$$\includegraphics[align=c,scale=1]{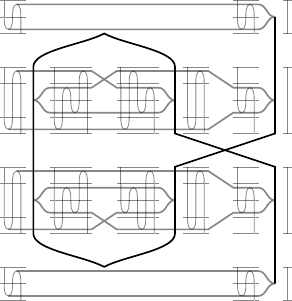}=\Id_{C_r}.$$
	
\end{lemma}

\begin{proof}
	
The proof is the following string diagram computation.
	
\begin{tabular}{llll}
	
	\includegraphics[align=c,scale=1]{adj/obj/extra/butterfly1.pdf}  & $=$ & & \\ \includegraphics[align=c,scale=1]{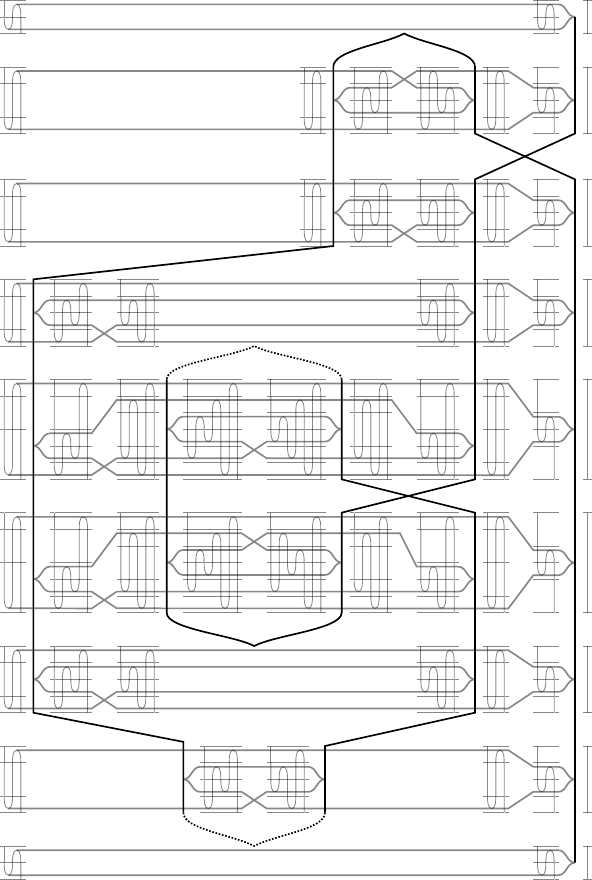} & $=$ & \includegraphics[align=c,scale=1]{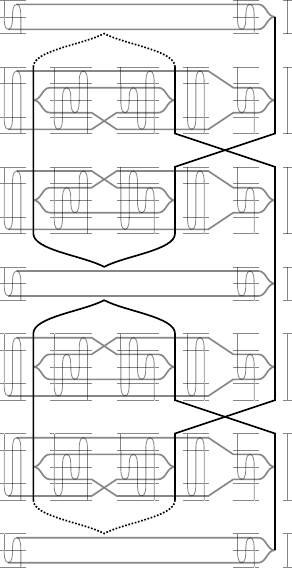} & $=$
	
\end{tabular}

\begin{tabular}{lllll}
	
	\includegraphics[align=c,scale=1]{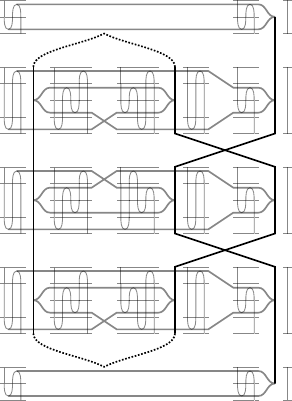} & $=$ & \includegraphics[align=c,scale=1]{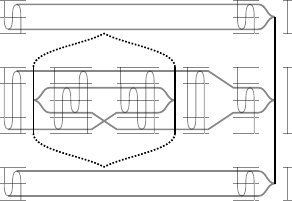} & $=$ &
	\includegraphics[align=c,scale=1]{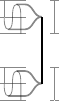}
	
\end{tabular} .

\end{proof}

\begin{lemma}\label{butterfly2}
	
	The $4$-morphisms $\SW^{-1}$ and $\SW_2$ (as defined above) satisfy the butterfly relation 	$$\includegraphics[align=c,scale=1]{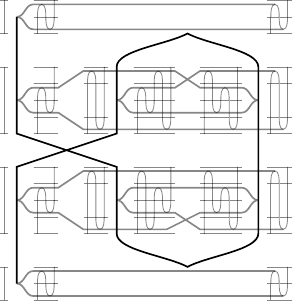}=\Id_{C_l^{-1}}.$$
	
\end{lemma}

\begin{proof}
	
The proof is the following string diagram computation.

\begin{center}
	\begin{tabular}{llll}
		
		\includegraphics[align=c,scale=1]{adj/obj/extra/butterfly2.pdf}  & $=$ & \includegraphics[align=c,scale=1]{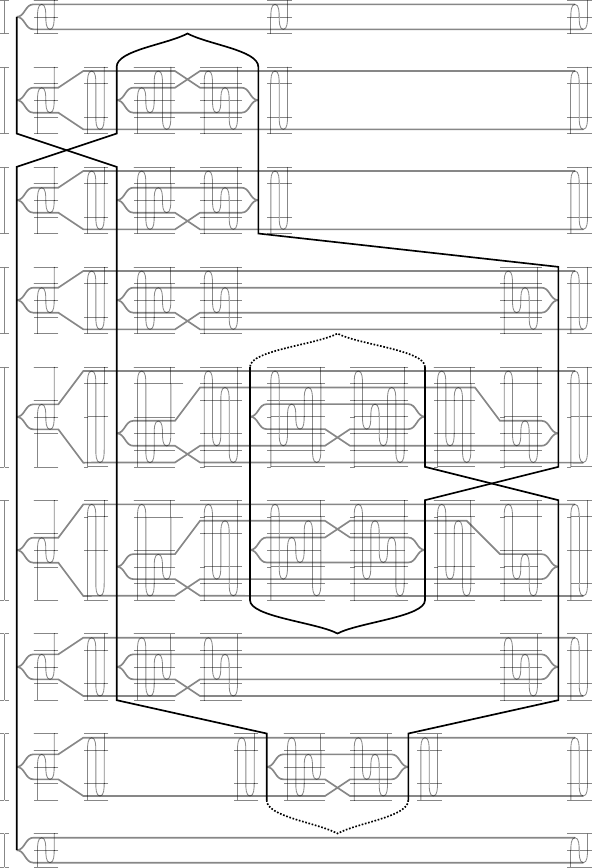} & $=$
		
	\end{tabular}
	
\end{center}

\begin{tabular}{lll}
	
	$=$ & \includegraphics[align=c,scale=1]{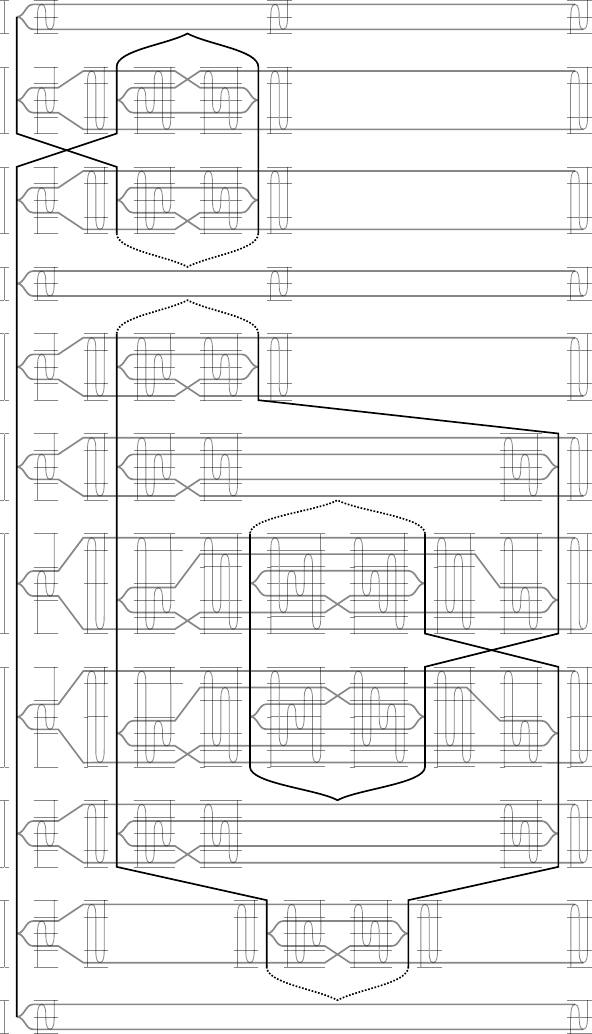} & $=$
	
\end{tabular}	
	
\begin{tabular}{lll}
	
	$=$ & \includegraphics[align=c,scale=1]{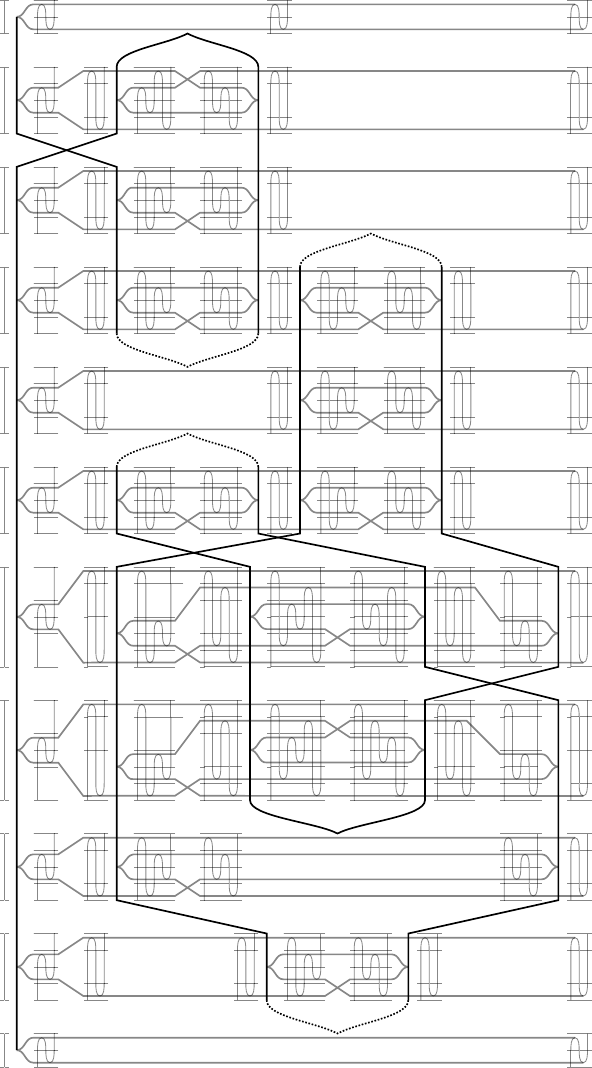} & $=$
	
\end{tabular}	
	
\begin{tabular}{lll}
	
	$=$ & \includegraphics[align=c,scale=1]{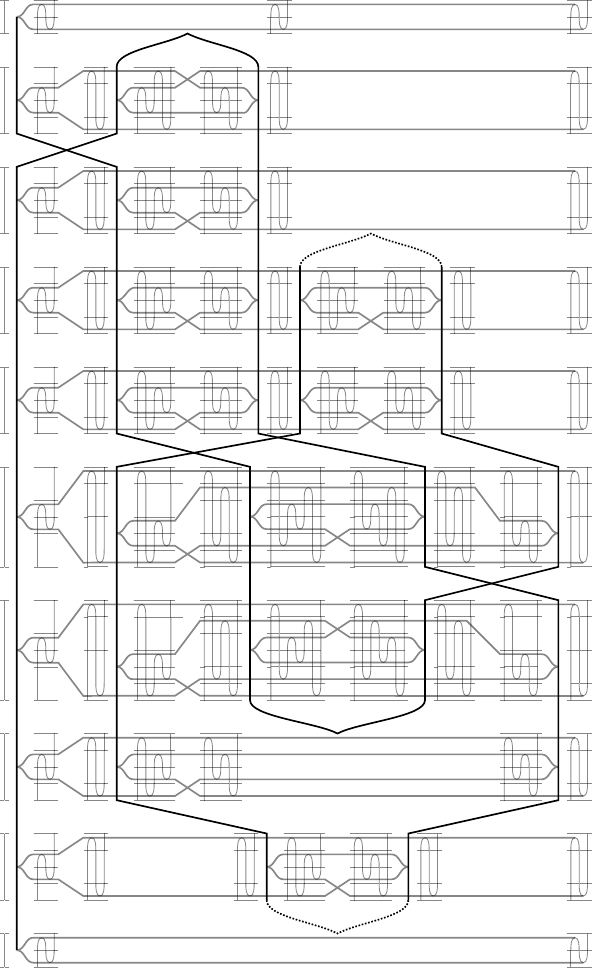} & $=$
	
\end{tabular}	
	
\begin{tabular}{lll}
	
	$=$ & \includegraphics[align=c,scale=1]{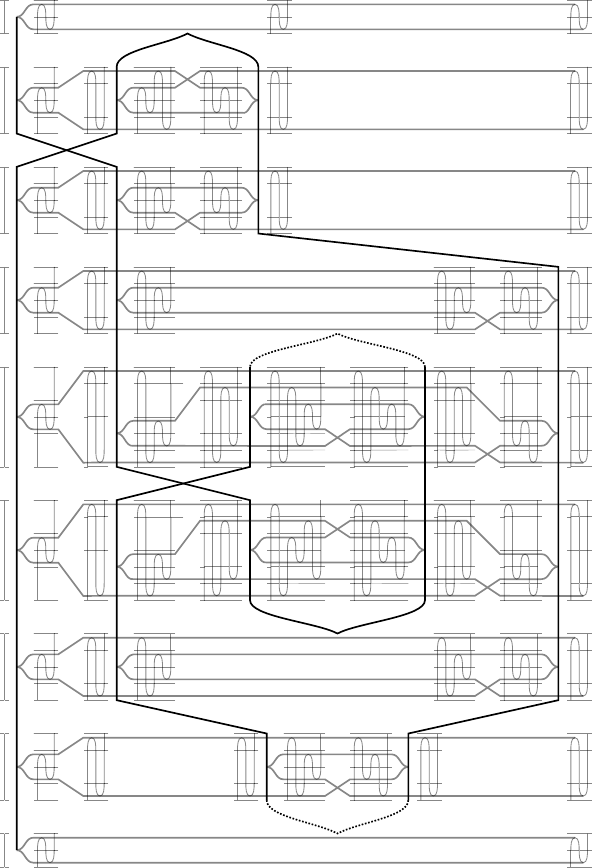} & $=$
	
\end{tabular}

\begin{tabular}{llllll}
	
	$=$ & \includegraphics[align=c,scale=1]{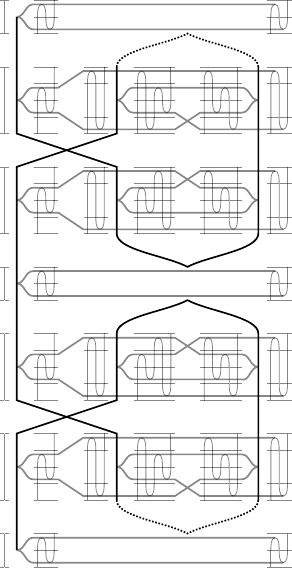} & $=$ & \includegraphics[align=c,scale=1]{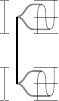} & $.$
	
\end{tabular}\end{proof}

Now we show that in fact the butterfly relations uniquely determine our choice of $\SW_2$.

\begin{definition}
	
Let $\Adj_{(4,1)}^{++}$ be the presentation obtained from $\Adj_{(4,1)}^+$ by adding generating $4$-cells $\SW_2$ and $\SW_2^{-1}$ and the two butterfly relations, as above, together with two relations witnessing the fact that $\SW_2$ and $\SW_2^{-1}$ are inverse to each other.
	
\end{definition}

\begin{lemma}

Let $\CC$ be a $4$-category and $F,G:F(\Adj_{(4,1)}^{++})\to \CC$ functors which agree on $\Adj_{(4,1)}^+$. Then $F=G$.

\end{lemma}

\begin{proof}
	
We need to show that $F(\SW_2)=G(\SW_2)$, which amounts to the fact that the butterfly relations already uniquely constrain the choice of $\SW_2$, once the choice of $\SW$ is fixed. This can be proved by a string diagram calculation which is very similar to the one in the proof of Lemma \ref{aswdef} below.	
\end{proof}

\begin{proposition}
	
The functor $F(\Adj_{(4,1)}^+)\to F(\Adj_{(4,1)}^{++})$ induced by the inclusion of presentations is an isomorphism.
	
\end{proposition}

\begin{proof}
	
The above constructions determine a functor $F(\Adj_{(4,1)}^{++})\to F(\Adj_{(4,1)}^+)$ which is the identity on $F(\Adj_{(4,1)}^+)$. So we only need to show that the composite $F(\Adj_{(4,1)}^{++})\to F(\Adj_{(4,1)}^+)\to F(\Adj_{(4,1)}^{++})$ is the identity. This follows from the above Lemma.
\end{proof}

\subsubsection{The swallowtail flip relations}

In addition to $\SW$ and $\SW_2$, we can define isomorphisms $\overline{\SW}$ and $\overline{\SW_2}$ between the inverses of the sources of $\SW$ and $\SW_2$ and the appropriate identity morphisms. These are related to $\SW,\SW_2$ by the swallowtail flip relations and they also satisfy two butterfly relations.

\begin{longtable}{lllll}
		
		$\overline{\SW}=$ & \includegraphics[align=c,scale=1]{adj/obj/extra/swallowtail2.pdf} $:$ \includegraphics[align=c,scale=1]{adj/obj/extra/swallowtail2_s.pdf} & $\xymatrix@1{\ar@1@<1ex>[r] & \ar@1@<1ex>[l]}$ & \includegraphics[align=c,scale=1]{adj/pres/id_u.pdf} $:$ \includegraphics[align=c,scale=1]{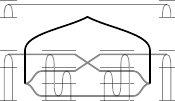} & $=\overline{\SW}^{-1}$ 
		
		\\
		
		\\
		
		$\overline{\SW_2}=$ &\includegraphics[align=c,scale=1]{adj/obj/extra/swallowtail4.pdf} $:$ \includegraphics[align=c,scale=1]{adj/obj/extra/swallowtail4_s.pdf} & $ \xymatrix@1{\ar@1@<1ex>[r] & \ar@1@<1ex>[l]}$ & \includegraphics[align=c,scale=1]{adj/pres/id_c.pdf} $:$ \includegraphics[align=c,scale=1]{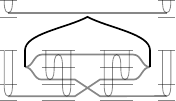} & $=\overline{\SW_2}^{-1}$

\end{longtable}	

Namely, we can take

\begin{center}
	
	$\overline{\SW}=$\includegraphics[align=c,scale=1]{adj/obj/extra/swallowtail2.pdf} $:=$ \includegraphics[align=c,scale=1]{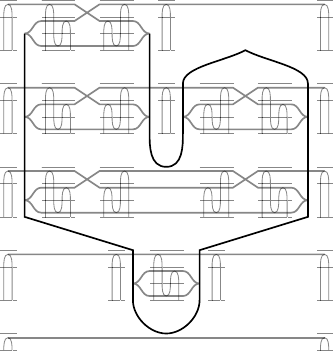}  
	
\end{center}

and 

\begin{center}
	
	$\overline{\SW_2}=$\includegraphics[align=c,scale=1]{adj/obj/extra/swallowtail4.pdf} $:=$ \includegraphics[align=c,scale=1]{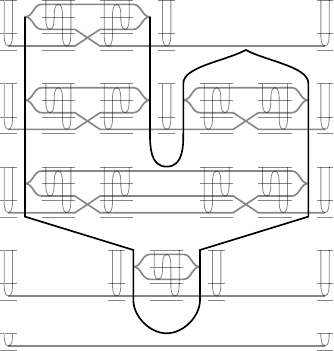}.  
	
\end{center}

\begin{remark}\label{flip}
	
We could also have defined $\overline{\SW}$ and $\overline{\SW_2}$ by bending the legs of $\SW$ and $\SW_2$ to the right, instead of the left. It is easy to show, using the snake relations, that the resulting composites would be equal to the ones we have written down.	
	
\end{remark}	

\begin{lemma}
	
The pairs $(\SW^{-1},\overline{\SW})$ and $(\SW_2^{-1},\overline{\SW_2})$ satisfy the two \textbf{swallowtail flip relations}: 

\begin{center}\begin{tabular}{ccc}
		
		\includegraphics[align=c,scale=1]{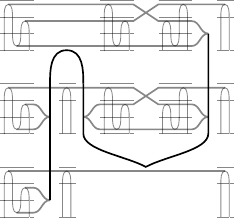}  & $=$ & \includegraphics[align=c,scale=1]{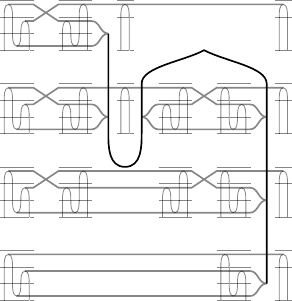} 
		
		\\
		
		\\
		
		\includegraphics[align=c,scale=1]{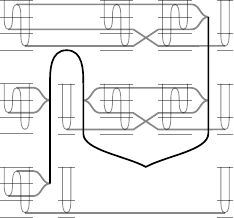}  & $=$ & \includegraphics[align=c,scale=1]{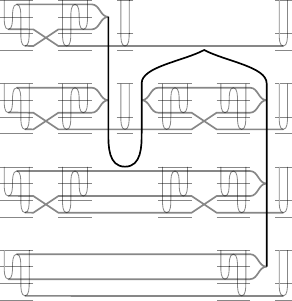}
		
	\end{tabular}\end{center}	
	
\end{lemma}

\begin{proof}
	
This is immediate from the definition of $\overline{\SW},\overline{\SW_2}$, using the snake relations.
\end{proof}		

\begin{lemma}
	
The $4$-morphisms $\overline{\SW}^{-1}$ and $\overline{\SW_2}$ as defined above satisfy the butterfly relations

\begin{center}\begin{tabular}{ccccc}

 \includegraphics[align=c,scale=1]{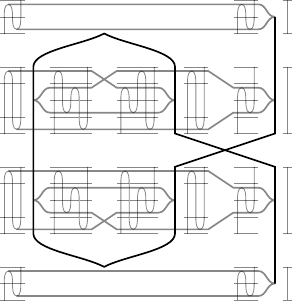} & $=\Id_{C_l}$ & and & \includegraphics[align=c,scale=1]{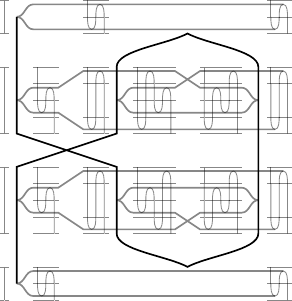} & $=\Id_{C_r^{-1}}.$

		\end{tabular}\end{center}

\end{lemma}

\begin{proof}
	
The proof of the first relation is the string diagram computation below. The proof of the second relation is a very similar computation.

\begin{longtable}{llll}
	
	\includegraphics[align=c,scale=1]{adj/obj/extra/butterfly3.pdf}  & $=$ &	\includegraphics[align=c,scale=1]{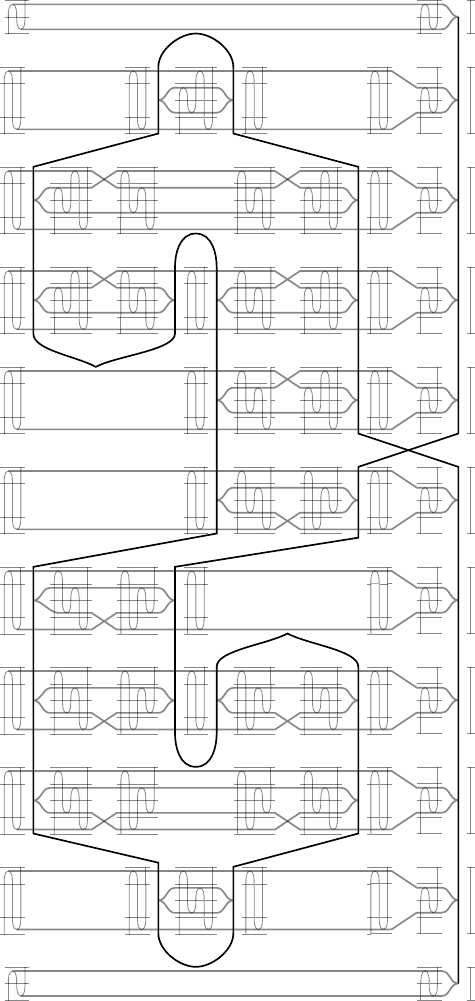} &  $=$
	
\end{longtable}
	
\begin{longtable}{llll}	
	
 \includegraphics[align=c,scale=1]{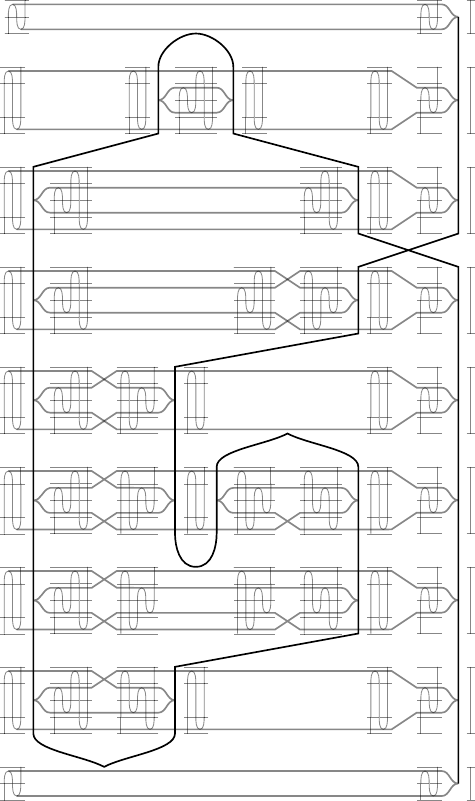} & $=$ \includegraphics[align=c,scale=1]{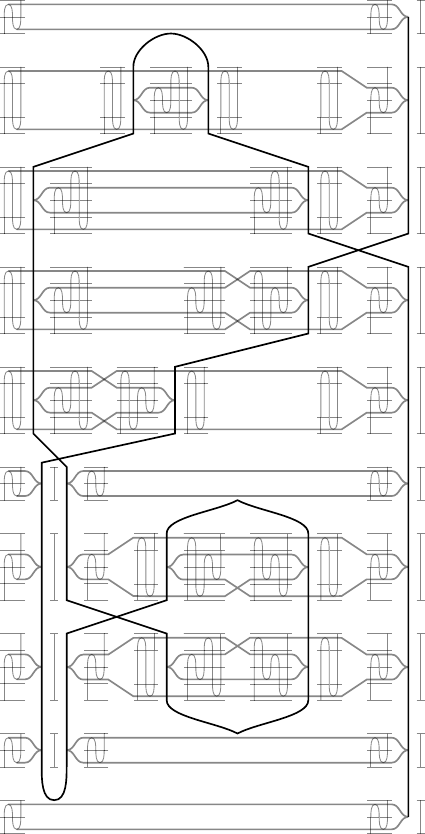} & $=$
 
\end{longtable}

\begin{longtable}{lllll}	
 
   \includegraphics[align=c,scale=1]{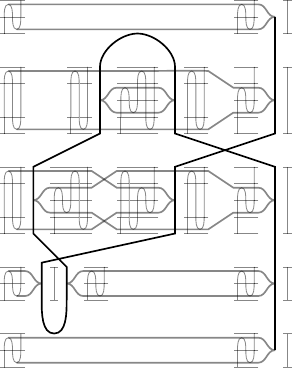} & $=$ & \includegraphics[align=c,scale=1]{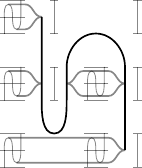} & $=$ & \includegraphics[align=c,scale=1]{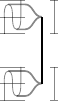} .

\end{longtable}\end{proof}

\begin{definition}
	
Let $\Adj_{(4,1)}^{+++}$ be the presentation obtained from $\Adj_{(4,1)}^{++}$ by adding generating $4$-cells $\overline{\SW}$, $\overline{\SW}^{-1}$, $\overline{\SW_2}$ and $\overline{\SW_2}^{-1}$ together with four relations witnessing the fact that these two pairs of cells are actually inverse two each other; the two butterfly relations, as above; plus the two swallowtail flip relations as above.

\end{definition}

\begin{lemma}
	
	Let $\CC$ be a $4$-category and $F,G:F(\Adj_{(4,1)}^{+++})\to \CC$ functors which agree on $\Adj_{(4,1)}^{++}$. Then $F=G$.
	
\end{lemma}

\begin{proof}
	
The swallowtail flip relations imply that the images of $\overline{\SW}$ and $\overline{\SW_2}$ under $F$ and $G$ must be defined from the images of $\SW^{-1}$ and $\SW_2^{-1}$ by the formulas we used above. Finally the images of $\overline{\SW}^{-1}$ and $\overline{\SW_2}^{-1}$ must be the inverses of the images of $\overline{\SW}$ and $\overline{\SW_2}$. This shows that	$F$ and $G$ are determined by their values on $\Adj_{(4,1)}^{++}$, so $F=G$.
\end{proof}

\begin{proposition}
	
The functor $F(\Adj_{(4,1)}^{++})\to F(\Adj_{(4,1)}^{+++})$ induced by the inclusion of presentations is an isomorphism.	
	
\end{proposition}

\begin{proof}
	
The constructions above determine a functor $F(\Adj_{(4,1)}^{+++})\to F(\Adj_{(4,1)}^{++})$ which is inverse to the inclusion, by the previous Lemma. 	\end{proof}	

\subsubsection{The final enlarged presentation}	

Finally, one can add some extra relations which already hold in $F(\Adj_{(4,1)}^{+++})$. There are two swallowtail flip relations that look exactly like the ones we wrote down, except we bend one leg to the right, instead of the left. These follow from the ones we wrote down, by Remark \ref{flip}. There are then four additional swallowtail flip relations, whose sources are inverse to the four already described. Finally there are four additional butterfly relations, whose sources are inverse to the four already described.

\begin{definition}
	
We define the presentation $\Adj_{(4,1)}^{max}$ by adding to $\Adj_{(4,1)}^{+++}$ the six swallowtail flip relations and four butterfly relations described above.

This presentation therefore includes a total of four inverse pairs of swallowtail $4$-morphisms (one of these pairs is already in $\Adj_{(4,1)}$), eight butterfly relations, eight swallowtail flip relations and eight snake relations (two of which are in $\Adj_{(4,1)}$). 	
	
\end{definition}

\begin{lemma}
	
We have $F(\Adj_{(4,1)}^{max})=F(\Adj_{(4,1)}^{+++})$	
	
\end{lemma}	

\begin{proof}
	
We have only added relations which already hold in $\Adj_{(4,1)}^{+++}$.	\end{proof}	

\begin{corollary}
	
The inclusion map $F(\Adj_{(4,1)})\to F(\Adj_{4,1}^{max})$ is an isomorphism, with an explicit inverse given by the constructions in this section. 
	
\end{corollary}	

In view of the above result, we will sometimes use generators and relations from $\Adj_{4,1}^{max}$ in our computations in the rest of this paper, as these can all be interpreted inside $F(\Adj_{(4,1)})$.

\section{Proof of the main Theorem}

We now give a proof of the Main Theorem, which we state here again, for convenience.

\begin{theorem}[Main Theorem]\label{main}
	
	Given a strict $4$-category $\CC$, the restriction map $$E_l:\Map(\Adj_{(4,1)},\CC)\to\Map^{L}(\theta^{(1)}, \CC)$$ is a trivial fibration of strict $4$-groupoids.
	
\end{theorem}

Theorem \ref{main} will be a consequence of the following Lemma.

\begin{lemma}\label{pullback}
 
Given a $4$-category $\CC$, the map $$\psi:\Map(\Adj_{(4,1)},\CC)\to\Map(\Adj_{(4,1)},\h_3(\CC))\times_{\Map^L(\theta^{(1)}, \h_3(\CC))}\Map^L(\theta^{(1)}, \CC)$$ induced by the square $$\xymatrix{\Map(\Adj_{(4,1)},\CC)\ar[d]\ar[r]^-{E_l} & \Map^L(\theta^{(1)}, \CC)\ar[d] \\ \Map(\Adj_{(4,1)},\h_3(\CC))\ar[r]_-{E_l} & \Map^L(\theta^{(1)}, \h_3(\CC))}$$ is a trivial fibration.

\end{lemma}

Given this Lemma, it is easy to give a proof of Theorem \ref{main}.

\begin{proof}[Proof of Theorem \ref{main}]
	
Note that the bottom map in the square from Lemma \ref{pullback} is the same as $$\Map(\Adj_{(3,1)},\h_3(\CC))\to\Map^L(\theta^{(1)}, \h_3(\CC))$$ which is a fibration by Theorem \ref{fibration} and a weak equivalence of $3$-groupoids by Proposition \ref{adj3}. Therefore, applying Lemma \ref{triv_fib}, the map $$\pi_2:\Map(\Adj_{(4,1)},\h_3(\CC))\times_{\Map^L(\theta^{(1)}, \h_3(\CC))}\Map^L(\theta^{(1)}, \CC)\to\Map^L(\theta^{(1)}, \CC)$$ is a trivial fibration of $4$-groupoids. By Lemmas \ref{pullback}, \ref{comp_weak_eq} and \ref{comp_fib}, the composite $\Map(\Adj_{(4,1)},\CC)\to\Map^L(\theta^{(1)},\CC)$ is a trivial fibration.\end{proof}

 So it is enough to prove Lemma \ref{pullback}. Since $\psi$ is a fibration (Corollary \ref{corfib}), it is enough to show that it has weakly contractible fibres. This is what we will do in the rest of this section.
 
\begin{lemma}\label{3groupoids}
  
The fibres of $\psi$ are $3$-groupoids.  
  
\end{lemma}

\begin{proof}

Note that a priori the fibres are $4$-groupoids, being fibres of a map between $4$-groupoids. So we have to show that the fibres contain only identitity $4$-morphisms. This follows from the fact that $\theta^{(1)}$ contains the $0$-skeleton of $\Adj_{(4,1)}$.\end{proof}

\subsection{The fibres are nonempty}\label{section_nonempty}

We show that $$\psi:\Map(\Adj_{(4,1)},\CC)\to\Map(\Adj_{(4,1)},\h_3(\CC))\times_{\Map^L(\theta^{(1)}, \h_3(\CC))}\Map^L(\theta^{(1)}, \CC)$$ is surjective on objects. In other words, we show its fibres are all nonempty.

Given $F:\Adj_{(4,1)}\to\h_3(\CC)$ we must lift it to a functor $\Adj_{(4,1)}\to\CC$.

\begin{definition}
 
We define $\Adj^-_{(4,1)}$ to be the presentation obtained from $\Adj_{(4,1)}$ by removing the two snake relations

\begin{longtable}{lll}

\includegraphics[align=c,scale=1]{adj/pres/snake_cl.pdf} $=\Id_{C_l}$ & and & \includegraphics[align=c,scale=1]{adj/pres/snake_cr.pdf} $=\Id_{C_r}.$

   \end{longtable}
 
\end{definition}
 
\begin{lemma}\label{nonempty}
 
The map $\psi$ is surjective on objects.
 
\end{lemma}
 
\begin{proof}
 
Consider $F:\Adj_{(4,1)}\to h_3(\CC)$. From this one can construct a functor $F^-:\Adj_{(4,1)}^-\to \CC$. In order to obtain a functor $F:\Adj_{(4,1)}\to\CC$ we apply Lemma \ref{adjeq} to the invertible $3$-morphisms $F(C_r)$ and $F(C_l)$, with inverses $F(C_r^{-1})$ and $F(C_l^{-1})$ and units $F(u_{C_r})$ and $F(u_{C_l})$. We obtain new counit $4$-morphisms satisfying the appropriate snake equations. Letting $F(c_{C_r})$ and $F(c_{C_l})$ equal these morphisms, we get a functor $F:\Adj_{(4,1)}\to\CC$ which lifts the original $F$.\end{proof}

\subsection{ The fibres are connected}\label{section_connected}

Now we prove that the fibres of $$\psi:\Map(\Adj_{(4,1)},\CC)\to\Map(\Adj_{(4,1)},\h_3(\CC))\times_{\Map^L(\theta^{(1)}, \h_3(\CC))}\Map^L(\theta^{(1)}, \CC)$$ are connected. 

\begin{lemma}\label{aswdef}
 
Given $F,G\in\Map(\Adj_{(4,1)},\CC)$ with $F=G$ on the $3$-skeleton of $\Adj_{(4,1)}$, there exists an isomorphism $\alpha:F\to G$ in $$\Map(\sk_3(\Adj_{(4,1)})\cup\{\SW,\SW^{-1}\},\CC)$$ which is the identity on $\sk_3(\Adj_{(4,1)})\setminus \{C_l^{-1}\}$. 
 
\end{lemma}

\begin{proof}
 
Since $F=G$ on $\{X,Y,l,r,\un,\co,C_l,C_l^{-1},C_r,C_r^{-1}\}$, the relation $\alpha_{\SW}$ takes the form \begin{longtable}{lll}\includegraphics[align=c,scale=1]{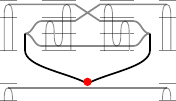} & $=$ & \includegraphics[align=c,scale=1]{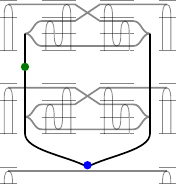} .\end{longtable} We define $\alpha_{C_l^{-1}}$ to be the $4$-morphism \begin{center}\includegraphics[align=c,scale=1]{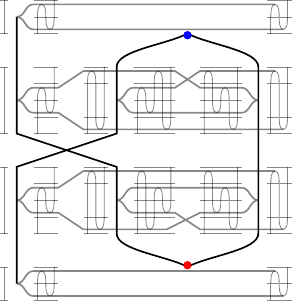}\end{center} and we check that the relation $\alpha_{\SW}$ holds as follows

\begin{center}
 \begin{longtable}{llll}
  
\includegraphics[align=c,scale=1]{adj/1morph/a_sw_proof_1.pdf}  & $=$ & \includegraphics[align=c,scale=1]{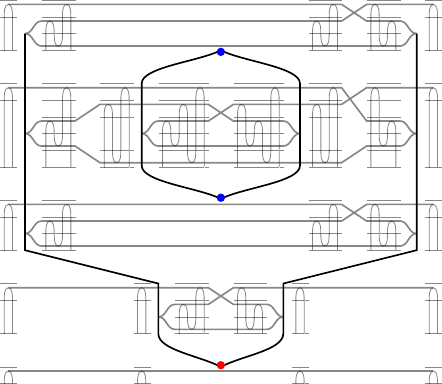}  & $=$ 

\\

\includegraphics[align=c,scale=1]{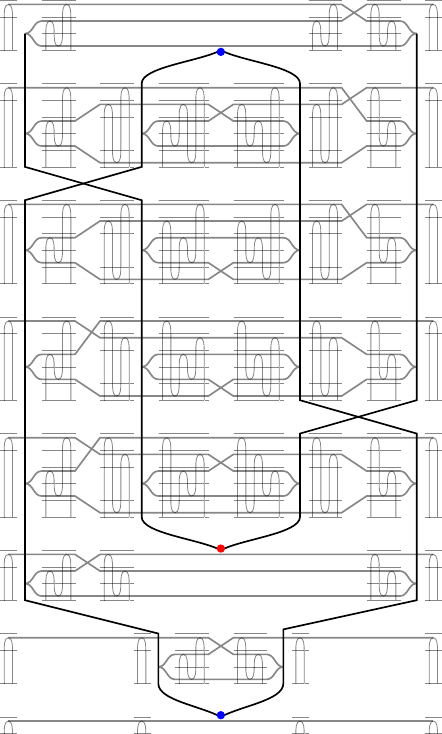}  & $=$ & \includegraphics[align=c,scale=1]{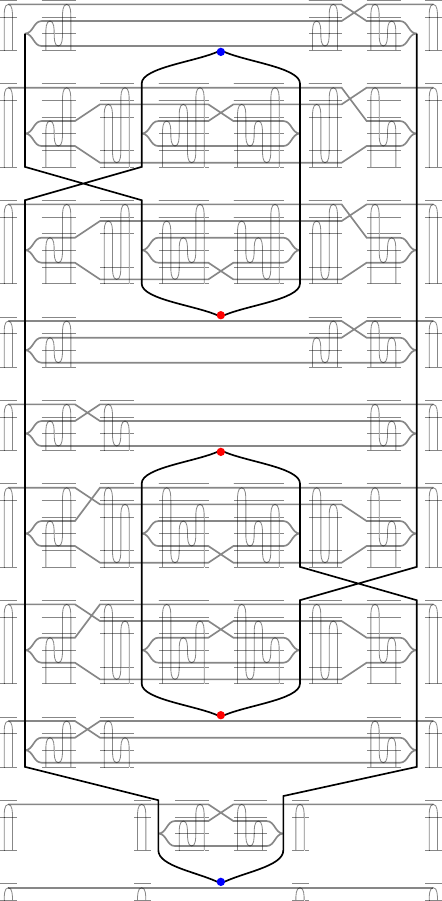}  & $=$ 

\\ \\

\includegraphics[align=c,scale=1]{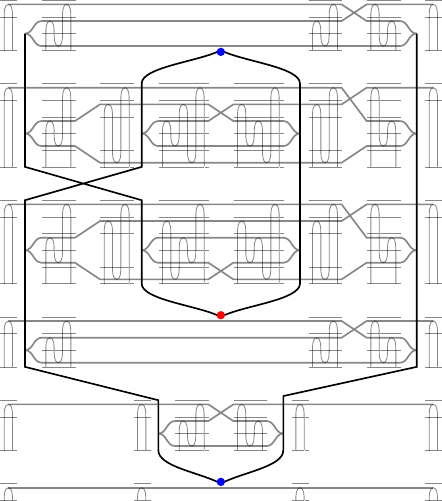}  & $.$

 \end{longtable}
 \end{center}
 
The relation $\alpha_{\SW^{-1}}$ follows from $\alpha_{\SW}$. \end{proof}

\begin{lemma}\label{afinaldef}
 
Given $F,G\in\Map(\Adj_{(4,1)},\CC)$ with $F=G$ on $$\sk_3(\Adj_{(4,1)})\cup\{\SW,\SW^{-1}\}$$ there exists an isomorphism $\alpha:F\to G$ in $\Map(\Adj_{(4,1)},\CC)$, which is the identity on $\sk_3(\Adj_{(4,1)})\setminus\{C_l,C_r^{-1}\}$.
 
\end{lemma}

\begin{proof}
 
We define \begin{longtable}{lllll}$\alpha_{C_l}:=$ & \includegraphics[align=c,scale=1]{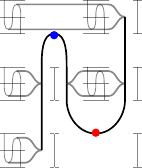} & and & $\alpha_{C_r^{-1}}:=$ & \includegraphics[align=c,scale=1]{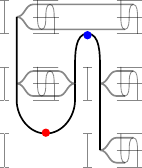}\end{longtable} so the relations  $\alpha_{u_{C_l}}$, $\alpha_{c_{C_l}}$ , $\alpha_{u_{C_r}}$ and $\alpha_{c_{C_r}}$ are satisfied. The relations $\alpha_{u_{C_l}^{-1}}$, $\alpha_{c_{C_l}^{-1}}$ , $\alpha_{u_{C_r}^{-1}}$ and $\alpha_{c_{C_r}^{-1}}$ follow from these ones. The relation $\alpha_{\SW}$ holds because $\alpha_{C_l^{-1}}$ and $\alpha_{C_r}$ are the identity and $F(\SW)=G(\SW)$. The relation $\alpha_{\SW^{-1}}$ follows from $\alpha_{\SW}$.\end{proof}

\begin{lemma}\label{connected_core}
 
Given $F,G\in\Map(\Adj_{(4,1)},\CC)$ such that $F=G$ on the $3$-skeleton, there exists $\alpha:F\to G$ in $\Map(\Adj_{(4,1)},\CC)$ such that $\alpha=\id_F$ on the $2$-skeleton.  
 
\end{lemma}

\begin{proof}By Lemma \ref{aswdef} there is an isomorphism $\alpha:F\to G$ in $$\Map(\sk_3(\Adj_{(4,1)})\cup\{\SW,\SW^{-1}\},\CC)$$ which is the identity on $\sk_3(\Adj_{(4,1)})\setminus \{C_l^{-1}\}$. Since the restriction map $$\Map(\Adj_{(4,1)},\CC)\to\Map(\sk_3(\Adj_{(4,1)})\cup\{\SW,\SW^{-1}\},\CC)$$ is a fibration (Theorem \ref{fibration}), one can extend this to an equivalence $\alpha:F\to F_1$ in $\Map(\Adj_{(4,1)},\CC)$. Now it is enough to find an equivalence $F_1\to G$ which is the identity on $\sk_2(\Adj_{(4,1)})$. Since $F_1=G$ on $\sk_3(\Adj_{(4,1)})\cup\{\SW,\SW^{-1}\}$, this can be done by Lemma \ref{afinaldef}.\end{proof}

\begin{lemma}\label{connected}
 
The fibres of $$\psi:\Map(\Adj_{(4,1)},\CC)\to\Map(\Adj_{(4,1)},\h_3(\CC))\times_{\Map^L(\theta^{(1)}, \h_3(\CC))}\Map^L(\theta^{(1)}, \CC)$$ are connected.
 
\end{lemma}

\begin{proof}

By Lemma \ref{nonempty}, the fibres are nonempty. So we just need to check that all objects in the same fibre are equivalent. Given $F,G\in\Map(\Adj_{(4,1)}, \CC)$ with $\psi(F)=\psi(G)$ we need to find an equivalence $\alpha:F\to G$ in $\Map(\Adj_{(4,1)},\CC)$ such that $\psi(\alpha)=\id_{\psi(F)}$. Equivalently, we are given $F,G\in\Map(\Adj_{(4,1)}, \CC)$ with $F=G$ on the $2$-skeleton of $\Adj_{(4,1)}$ and $F(f)\simeq G(f)$ for every $3$-cell $f\in\Adj_{(4,1)}$. We need to find $\alpha:F\to G$ in $\Map(\Adj_{(4,1)}, \CC)$ such that $\alpha=\id_F$ over the $1$-skeleton and $\alpha_x\simeq\id_{F(x)}$ for every $2$-cell $x$ in $\Adj_{(4,1)}$. 

We use the equivalences $F(f)\simeq G(f)$ to build $\alpha:F\to G$ over the $3$-skeleton of $\Adj_{(4,1)}$, with $\alpha=\id_F$ over the $2$-skeleton. Because $\Map(\Adj_{(4,1)},\CC)\to\Map(\sk_3(\Adj_{(4,1)}),\CC)$ is a fibration (Theorem \ref{fibration}), we can lift this to $\tilde{\alpha}:F\to\tilde{G}$ in $\Map(\Adj_{(4,1)},\CC)$. On the $3$-skeleton of $\Adj_{(4,1)}$ we have $\tilde{G}=G$ and $\tilde{\alpha}=\alpha$, so $\tilde{\alpha}=\id_F$ on the $2$-skeleton. 

Now it's enough to find an equivalence $\beta:\tilde{G}\to G$ in $\Map(\Adj_{(4,1)},\CC)$ which is the identity on the $2$-skeleton, as its composite with $\tilde{\alpha}$ will provide an equivalence $F\to G$ over $\Adj_{(4,1)}$ which is the identity over the $2$-skeleton, which is even more than we need. We can find such $\beta$ by applying Lemma \ref{connected_core}.\end{proof}

\subsection{The fibres are $1$-connected}\label{section_1connected}

Now we prove that the fibres of $$\psi:\Map(\Adj_{(4,1)},\CC)\to\Map(\Adj_{(4,1)},\h_3(\CC))\times_{\Map^L(\theta^{(1)}, \h_3(\CC))}\Map^L(\theta^{(1)}, \CC)$$ are $1$-connected.

\begin{lemma}\label{mcldef}
 
Given a $1$-morphism $\alpha:F\to F$ in $\Map(\Adj_{(4,1)},\CC)$ such that $\alpha$ is the identity on the $2$-skeleton, there exists a $2$-morphism $m:\alpha\to\Id_F$ in $\Map(\sk_2(\Adj_{(4,1)})\cup\{C_l,C_l^{-1}\},\CC)$ which is the identity on $\sk_2(\Adj_{(4,1)})\setminus\{\un\}$.  
 
\end{lemma}

\begin{proof}
 
Denote $\alpha_{C_l}:F(C_l)\to F(C_l)$ by \begin{center}\includegraphics[align=c,scale=1.2]{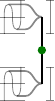}.\end{center} We want to define a $4$-morphism $m_u:\Id_u\to\Id_u$, which we denote by \begin{center}\includegraphics[align=c,scale=1.2]{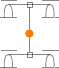},\end{center} such that the relations $m_{C_l}$ and $m_{C_l^{-1}}$ hold. The relation $m_{C_l}$ looks like \begin{longtable}{lll}\includegraphics[align=c,scale=1]{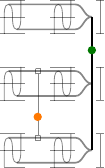} & $=$ & \includegraphics[align=c,scale=1]{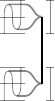}\end{longtable} and the relation $m_{C_l^{-1}}$ follows from $m_{C_l}$. We have 

\begin{longtable}{llll}
 
\includegraphics[align=c,scale=1]{adj/2morph/m_cl_proof_1.pdf} & $=$ & \includegraphics[align=c,scale=1]{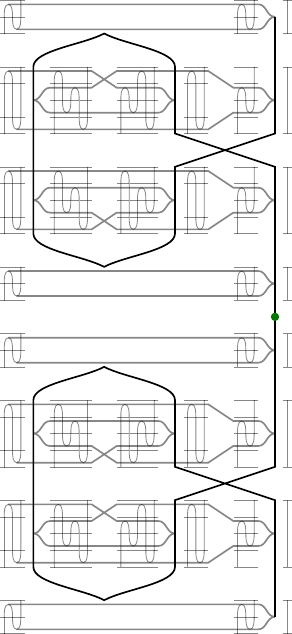} & $=$

\\

\\

\includegraphics[align=c,scale=1]{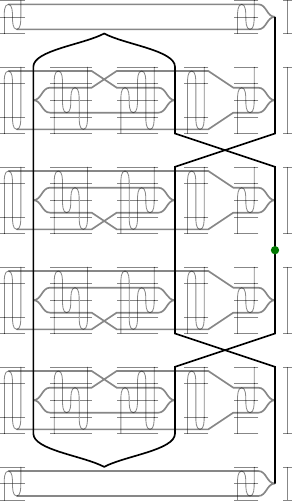} & $=$ & \includegraphics[align=c,scale=1]{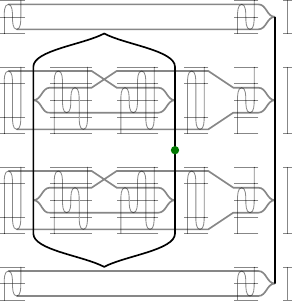} & 
 
\end{longtable}
 
so we just define $m_u$ to be the inverse of \begin{center}\includegraphics[align=c,scale=1]{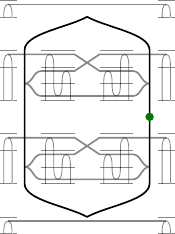}.\end{center}\end{proof}

\begin{lemma}\label{mcrdef}
 
Given a $1$-morphism $\alpha:F\to F$ in $\Map(\Adj_{(4,1)},\CC)$ such that $\alpha$ is the identity on $\sk_3(\Adj_{(4,1)})\setminus\{C_r,C_r^{-1}\}$, we have $\alpha=\Id$ in $\Map(\Adj_{(4,1)},\CC)$,
 
\end{lemma}

\begin{proof}

Consider the relation \begin{longtable}{llll}$\alpha_{SW}:$ & \includegraphics[align=c,scale=1]{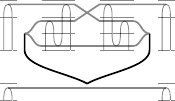} & $=$ & \includegraphics[align=c,scale=1]{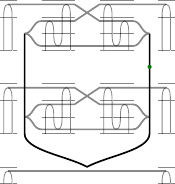}.\end{longtable} We use this to show $\alpha_{C_r}=\Id$ as follows.

\begin{center}
 
\begin{tabular}{llllll}
 
\includegraphics[align=c,scale=1]{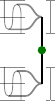} & $=$ & \includegraphics[align=c,scale=1]{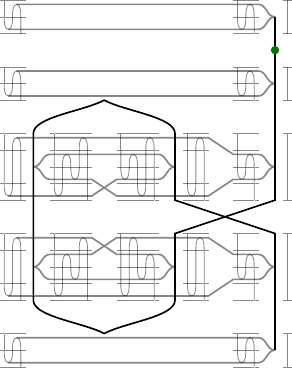} & $=$ & \includegraphics[align=c,scale=1]{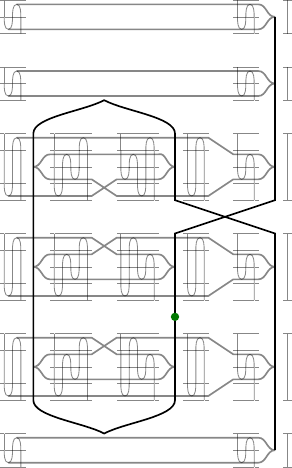} & $=$

\end{tabular}

\begin{tabular}{llll}

 $=$ & \includegraphics[align=c,scale=1]{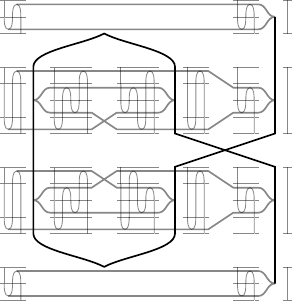} & $=$ & \includegraphics[align=c,scale=1]{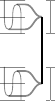} . 
 
\end{tabular}
 
\end{center}

It follows that $\alpha_{C_r^{-1}}$ is the identity, since $C_r\dashv C_r^{-1}$.\end{proof}

\begin{lemma}\label{1connected_core}

Given $F\in\Map(\Adj_{(4,1)},\CC)$ and an equivalence $\alpha:F\to F$ which is the identity on the $2$-skeleton, there exists a $2$-equivalence $m:\alpha\to \Id$ in $\Map(\Adj_{(4,1)},\CC)$ which is the identity on the $1$-skeleton.

\end{lemma}

\begin{proof}By Lemma \ref{mcldef}, there exists a $2$-morphism $m:\alpha\to\Id_F$ in $\Map(\sk_2(\Adj_{(4,1)})\cup\{C_l,C_l^{-1}\},\CC)$ which is the identity on $\sk_1(\Adj_{(4,1)})\cup\{\co\}$. Since $$\Map(\Adj_{(4,1)},\CC)\to\Map(\sk_2(\Adj_{(4,1)})\cup\{C_l,C_l^{-1}\},\CC)$$ is a fibration (Theorem \ref{fibration}), we can extend this to a $2$-morphism $m:\alpha\to\alpha_1$ in $\Map(\Adj_{(4,1)},\CC)$. Now $\alpha_1=\Id_F$ on $\sk_2(\Adj_{(4,1)})\cup\{C_l,C_l^{-1}\}$ so, by Lemma \ref{mcrdef}, we have $\alpha_1=\Id$ in $\Map(\Adj_{(4,1)},\CC)$.\end{proof}

\begin{lemma}\label{1connected}

The fibres of $\psi$ are $1$-connected.
 
\end{lemma}

\begin{proof}
 
By Lemma \ref{connected}, the fibres are connected, so we just need to show that $\pi_1(\psi^{-1}(\psi(F)),F)=*$ for any $F\in\Map(\Adj_{(4,1)},\CC)$. By an argument which is completely analogous to the one given in the proof of Lemma \ref{connected}, this reduces to Lemma \ref{1connected_core}. \end{proof}

\subsection{The fibres are $2$-connected}\label{section_2connected}

Now we prove that the fibres of $$\psi:\Map(\Adj_{(4,1)},\CC)\to\Map(\Adj_{(4,1)},\h_3(\CC))\times_{\Map^L(\theta^{(1)}, \h_3(\CC))}\Map^L(\theta^{(1)}, \CC)$$ are $2$-connected.

\begin{lemma}\label{aaudef}
 
Given a $2$-morphism $m:\Id_F\to \Id_F$ in $\Map(\Adj_{(4,1)},\CC)$ such that $m$ is the identity on the $1$-skeleton, there exists a $3$-morphism $\A:m\to\Id^{(2)}_F$ in $\Map(\sk_1(\Adj_{(4,1)})\cup\{\un\},\CC)$ which is the identity on $\{X,Y,l\}$.  
 
\end{lemma}

\begin{proof}
 
We need to define a $4$-morphism \begin{center}$\A_r:=$ \includegraphics[align=c,scale=1]{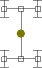} $:\id^{(2)}_{F(r)}\to \id^{(2)}_{F(r)}$\end{center} such that the relation \begin{center}$\A_u:$ \includegraphics[align=c,scale=1]{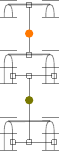} $=\id^{(2)}_{F(\un)}$ \end{center} holds. Here we denote \begin{center}$m_u:=$ \includegraphics[align=c,scale=1]{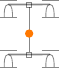}.\end{center} 

We have 
 
\begin{longtable}{llllll}
 
\includegraphics[align=c,scale=1]{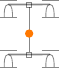} $=$   \\ \\

 $=$ \includegraphics[align=c,scale=1]{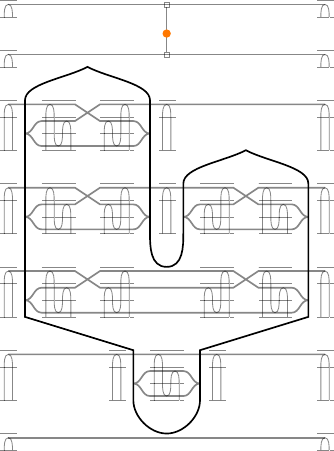} & $=$ &

\includegraphics[align=c,scale=1]{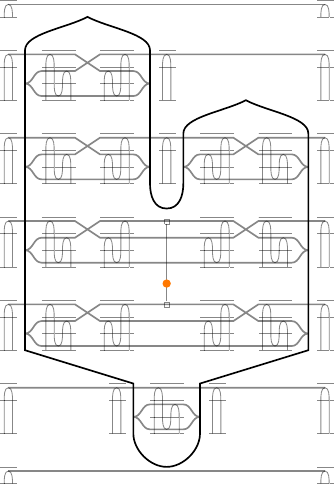} & $=$ &

\includegraphics[align=c,scale=1]{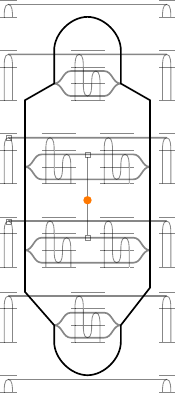} & .
 
\end{longtable}

so we just define $\A_r$ to be the inverse of \begin{center}\includegraphics[align=c,scale=1]{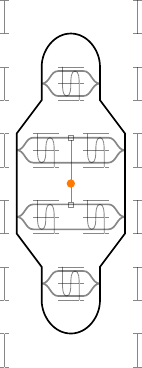}.\end{center} \end{proof}

\begin{lemma}\label{aacdef}
 
Given a $2$-morphism $m:\Id_F\to \Id_F$ in $\Map(\Adj_{(4,1)},\CC)$ such that $m$ is the identity on $\sk_2(\Adj_{(4,1)})\setminus\{\co\}$, we have $m=\Id$ in $\Map(\Adj_{(4,1)},\CC)$.
 
\end{lemma}

\begin{proof}
 
Denote \begin{center}$m_c:=$ \includegraphics[align=c,scale=1]{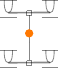}\end{center} and consider the relation \begin{center}$m_{C_l}:$ \includegraphics[align=c,scale=1]{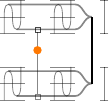} $=$ \includegraphics[align=c,scale=1]{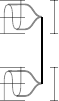}.\end{center} Then we have

\begin{longtable}{lllllllll}
 
\includegraphics[align=c,scale=1]{adj/3morph/a_c_proof_1.pdf} & $=$ &

\includegraphics[align=c,scale=1]{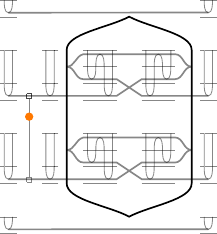} & $=$ &

\includegraphics[align=c,scale=1]{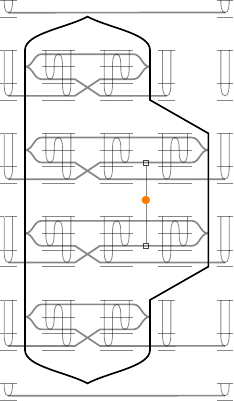} & $=$ &

\includegraphics[align=c,scale=1]{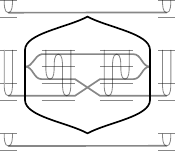} & $=$ &

\includegraphics[align=c,scale=1]{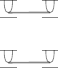} .  
 
\end{longtable}\end{proof}

\begin{lemma}\label{2connected_core}
 
Given $F:\Adj_{(4,1)}\to\CC$ and a $2$-morphism $m:\Id_F\to\Id_F$ in $\Map(\Adj_{(4,1)},\CC)$ such that $m$ is the identity on the $1$-skeleton, there exists a $3$-morphism $\A:m\to\Id^{(2)}_F$ in $\Map(\Adj_{(4,1)},\CC)$, which is the identity on $\{X,Y,l\}$. 
 
\end{lemma}

\begin{proof}By Lemma \ref{aaudef} there exists a $3$-morphism $\A:m\to\Id^{(2)}_F$ in $\Map(\sk_1(\Adj_{(4,1)})\cup\{\un\},\CC)$ which is the identity on $\{X,Y,l\}$. Since the restriction map $$\Map(\Adj_{(4,1)},\CC)\to\Map(\sk_1(\Adj_{(4,1)})\cup\{\un\},\CC)$$ is a fibration (Theorem \ref{fibration}), we can extend this to $\A:m\to m_1$ in $\Map(\Adj_{(4,1)},\CC)$. Now $m_1=\Id^{(2)}_F$ on $\sk_1(\Adj_{(4,1)})\cup\{\un\}$, so by Lemma \ref{aacdef} we have $m=\Id$ in $\Map(\Adj_{(4,1)},\CC)$.   \end{proof}

\begin{lemma}\label{2connected}

The fibres of $\psi$ are $2$-connected.
 
\end{lemma}

\begin{proof}
 
By Lemma \ref{1connected}, the fibres are $1$-connected, so we just need to show that $\pi_2(\psi^{-1}(\psi(F)),F)=*$ for any $F\in\Map(\Adj_{(4,1)},\CC)$. By an argument analogous to the one given in the proof of Lemma \ref{connected}, this reduces to Lemma \ref{2connected_core}.\end{proof}

\subsection{The fibres are $3$-connected}\label{section_3connected}

Now we prove that the fibres of $$\psi:\Map(\Adj_{(4,1)},\CC)\to\Map(\Adj_{(4,1)},\h_3(\CC))\times_{\Map^L(\theta^{(1)}, \h_3(\CC))}\Map^L(\theta^{(1)}, \CC)$$ are $3$-connected.

\begin{lemma}\label{3connected_core}
 
Given a $3$-morphism $\A:\Id^{(2)}_F\to \Id^{(2)}_F$ in $\Map(\Adj_{(4,1)},\CC)$ such that $\A$ is the identity on $\{X,Y,l\}$, we have $\A=\Id$ in $\Map(\Adj_{(4,1)},\CC)$.
 
\end{lemma}

\begin{proof}
 
All we have to do is show that \begin{center}$\A_r=$ \includegraphics[align=c,scale=1]{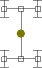} $:\Id^{(2)}_{F(r)}\to\Id^{(2)}_{F(r)}$\end{center} is the identity. Consider the relation \begin{center}$\A_{\un}:$ \includegraphics[align=c,scale=1]{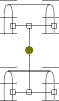} $=$ \includegraphics[align=c,scale=1]{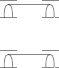}.\end{center} Using this, we have 

\begin{center}
 
\begin{longtable}{lllllllllll}
 
\includegraphics[align=c,scale=1]{adj/4morph/f_r_proof_1.pdf} & $=$ &

\includegraphics[align=c,scale=1]{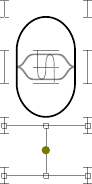} & $=$ &

\includegraphics[align=c,scale=1]{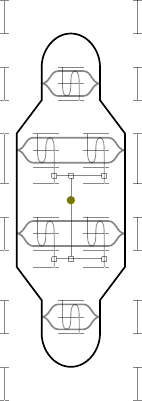} & $=$ &

\includegraphics[align=c,scale=1]{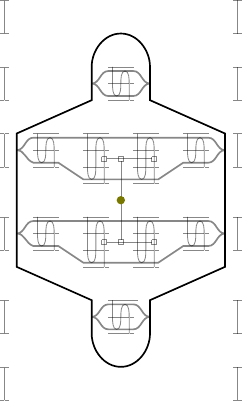} & $=$ &

\includegraphics[align=c,scale=1]{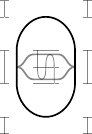} & $=$ &

\includegraphics[align=c,scale=1]{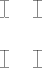} .
 
\end{longtable}
 
\end{center}\end{proof}

\begin{lemma}\label{3connected}

The fibres of $\psi$ are $3$-connected.
 
\end{lemma}

\begin{proof}
 
This follows from the Lemmas \ref{2connected} and \ref{3connected_core}.\end{proof}

\begin{proposition}\label{contractible}
 
The fibres of $\psi$ are weakly contractible. 
 
\end{proposition}

\begin{proof}
 
The fibres are $3$-groupoids (Lemma \ref{3groupoids}) and we have shown they are $3$-connected (Lemma \ref{3connected}). So by Corollary \ref{nconnected} the fibres are weakly contractible. \end{proof}

\begin{proof}[Proof of Lemma \ref{pullback}]

The map $$\psi:\Map(\Adj_{(4,1)},\CC)\to\Map(\Adj_{(4,1)},\h_3(\CC))\times_{\Map^L(\theta^{(1)}, \h_3(\CC))}\Map^L(\theta^{(1)}, \CC)$$ is a fibration by Corollary \ref{corfib} and its fibres are  weakly contractible by Proposition \ref{contractible}. Therefore, by Proposition \ref{fibre}, the map $\psi$ is a trivial fibration.\end{proof}

\section{Other versions of the main theorem}\label{more_main}

\begin{definition}
We define $\Unit$ to be the subcomputad of $\Adj_{(4,1)}$ containing the cells $\{X,Y,l,r,\un\}$. 
\end{definition}

\begin{definition}

We define $\Cusp$ to be the subcomputad of $\Adj_{(4,1)}$ containing $\Unit$ plus $\{\co,C_r\}$.
 
\end{definition}

\begin{definition}

We define $\Sw$ to be the subcomputad of $\Adj_{(4,1)}$ containing $\Cusp$ plus $\{C_l^{-1},\SW\}$.

\end{definition}

\begin{definition}
 We define $\Map^{\adj}(\Unit,\CC)$ to be the full $4$-subgroupoid of $\Map(\Unit,\CC)$ whose objects are those $F$ such that $(F(l),F(r),F(\un))$ can be extended to an adjunction in $\CC$.
\end{definition}

\begin{definition}

Define $\Map^{\adj}(\Cusp,\CC)$ to be the full $4$-subgroupoid of $\Map(\Cusp,\CC)$ whose objects are those $F$ such that $F(C_r)$ is an equivalence and $(F(l),F(r),F(u),F(c))$ is an adjunction.

\end{definition}

\begin{definition}

Define $\Map^{\adj}(\Sw,\CC)$ to be the full $4$-subgroupoid of $\Map(\Sw,\CC)$ whose objects are those $F$ such that $F(C_r)$, $F(C_l^{-1})$ and $F(\SW)$ are equivalences.
 
\end{definition}

\begin{theorem}
 
Let $\CC$ be a strict $4$-category. All restriction maps below are trivial fibrations of $4$-groupoids.

\[\begin{tikzcd}[column sep=tiny]
	{\Map(\Adj_{(4,1)},\CC)} & {\Map^{\adj}(\Sw,\CC)} & {\Map^{\adj}(\Cusp,\CC)} & {\Map^{\adj}(\Unit,\CC)} & {\Map^{L}(\theta^{(1)},\CC)}
	\arrow[from=1-1, to=1-2]
	\arrow[from=1-2, to=1-3]
	\arrow[from=1-3, to=1-4]
	\arrow[from=1-4, to=1-5]
\end{tikzcd}\]
 
\end{theorem}

In this section we will prove that the first map $$\phi:\Map^{\adj}(\Unit,\CC)\to\Map^L(\theta^{(1)},\CC)$$ is a trivial fibration. The proofs that the second and third maps are trivial fibrations are very similar and we ommit them. Finally the fact that the first map is a trivial fibration follows from the following simple argument.

\begin{proposition}
 
The map $E_{\Sw}:\Map(\Adj_{(4,1)},\CC)\to\Map^{\adj}(\Sw,\CC)$ is a trivial fibration. 
 
\end{proposition}

\begin{proof}

Consider \[\begin{tikzcd}[column sep=small]
	{\Map(\Adj_{(4,1)},\CC)} && {\Map^{\adj}(\Sw,\CC)} \\
	& {\Map^L(\theta^{(1)},\CC)}
	\arrow["E_{\Sw}",two heads, from=1-1, to=1-3]
	\arrow["E_l"', two heads, from=1-1, to=2-2]
	\arrow["\psi",from=1-3, to=2-2]
\end{tikzcd}.\]

We know $E_l$ is a trivial fibration by Theorem \ref{main}. The map $E_{\Sw}$ is a fibration because $\Map(\Adj_{(4,1)},\CC)\to\Map(\Sw,\CC)$ is a fibration by Theorem \ref{fibration}. Now $\psi$ is a weak equivalence because it's a composite of weak equivalences (Lemma \ref{comp_weak_eq}). Since weak equivalences satisfy the $2$ out of $3$ property (follows from Lemma \ref{htpygrps}) we conclude that $E_{\Sw}$ is also a weak equivalence.\end{proof}

In the rest of this section we prove that $\phi$ is a trivial fibration. Consider \[\begin{tikzcd}[column sep=small]
	{\Map(\Adj_{(4,1)},\CC)} && {\Map^{\adj}(\Unit,\CC)} \\
	& {\Map^L(\theta^{(1)},\CC)}
	\arrow["E_{\Unit}",two heads, from=1-1, to=1-3]
	\arrow["E_l"', two heads, from=1-1, to=2-2]
	\arrow["\phi",from=1-3, to=2-2]
\end{tikzcd}.\]

\begin{lemma}\label{k1}

Given $F\in\Map^{\adj}(\Unit,\CC)$ and $F\to G$ in $\Map(\Unit,\CC)$ we have $G\in\Map^{\adj}(\Unit,\CC)$.
 
\end{lemma}

\begin{proof}

Let $\PP$ be the subcomputad of $\Adj_{(4,1)}$ obtained by removing the swallowtail $3$-cells (and all relations involving them) and the snake relations between the cusp cancellation and creation $4$-morphisms. In other words, $\PP$ consists of $(l,r,\un,\co)$ together with the four cusp $3$-cells $(C_l,C_l^{-1},C_r,C_r^{-1})$, and all $4$-cells and relations witnessing the fact $C_r^{-1}$ is a weak inverse of $C_r$ and $C_l^{-1}$ is a weak inverse of $C_l$. Now $F\in\Map^{\adj}(\Unit,\CC)$ implies $F$ can be extended to $\PP$. Since $\Map(\PP,\CC)\to\Map(\Unit,\CC)$ is a fibration (Theorem \ref{fibration}), we can extend $F\to G$ to $\PP$. This gives an extension of $G$ to $\PP$, which shows that $G\in\Map^{\adj}(\Unit,\CC)$.\end{proof}

\begin{lemma}
 
The map $\phi$ is a fibration. 
 
\end{lemma}
\begin{proof}

In the diagram below the dotted arrow exists because $\Map(\Unit,\CC)\to\Map(\theta^{(1)},\CC)$ is a fibration (by Theorem \ref{fibration}). The dashed arrow exists because $\Map^{\adj}(\Unit,\CC)\hookrightarrow \Map(\Unit,\CC)$ is a full $4$-subgroupoid (together with Lemma \ref{k1} in the $k=1$ case). \[\begin{tikzcd}
	{\theta^{(k-1)}} & {\Map^{\adj}(\Unit,\CC)} & {\Map(\Unit,\CC)} \\
	{\theta^{(k)}} & {\Map^L(\theta^{(1)},\CC)} & {\Map(\theta^{(1)},\CC)}
	\arrow[from=1-1, to=1-2]
	\arrow[from=1-1, to=2-1]
	\arrow[hook, from=1-2, to=1-3]
	\arrow["\phi"'{pos=0.3}, from=1-2, to=2-2]
	\arrow[two heads, from=1-3, to=2-3]
	\arrow[dashed, from=2-1, to=1-2]
	\arrow[dotted, from=2-1, to=1-3]
	\arrow[from=2-1, to=2-2]
	\arrow[hook, from=2-2, to=2-3]
\end{tikzcd}\] \end{proof}

\begin{lemma} The map $\phi$ is surjective on objects. 
\end{lemma}
\begin{proof}An object in $\Map^L(\theta^{(1)},\CC)$ is by definition a $1$-morphism in $\CC$ which admits a right adjoint.\end{proof}

\begin{lemma}\label{connected_again}

The fibres of $\phi$ are connected.
 
\end{lemma}

\begin{proof} Given $F,G\in\Map^{\adj}(\Unit,\CC)$ such that $\phi(F)$=$\phi(G)$ we must find an equivalence $\alpha:F\to G$ in $\Map^{\adj}(\Unit,\CC)$ with $\alpha=\Id$ on $\{X,Y,l\}$. From $F,G$ we obtain $F,G\in\Map(\Adj_{(2,1)},\h_2(\CC))$. The fibres of $\Map(\Adj_{(2,1)},\h_2(\CC))\to\Map^L(\theta^{(1)},\h_2(\CC))$ are connected, by \cite[Lemma 3.10]{adj3}, so we have an equivalence $\alpha:F\to G$ in $\Map(\Adj_{(2,1)},\h_2(\CC))$ with $\alpha=\Id$ on $\{X,Y,l\}$. This consists of a $2$-morphism $\alpha_r$ and two equations $\alpha_u$ and $\alpha_c$. The equation $\alpha_u$ in $\h_2(\CC)$ can be lifted to a $3$-equivalence in $\CC$, so we obtain the desired $\alpha:F\to G$ in $\Map^{\adj}(\Unit,\CC)$. \end{proof}

\begin{lemma}\label{lrusurj}

The map $E_{\Unit}:\Map(\Adj_{(4,1)},\CC)\to\Map^{\adj}(\Unit,\CC)$ is surjective on objects.
 
\end{lemma}

\begin{proof}Given $G\in\Map^{\adj}(\Unit,\CC)$, there exists $\tilde{F}\in\Map(\Adj_{(4,1)},\CC)$ such that $E_l(\tilde{F})=\phi(G)$, because $E_{l}:\Map(\Adj_{(4,1)},\CC)\to\Map^L(\theta^{(1)},\CC)$ is surjective on objects (Theorem \ref{main}). Letting $F=E_{\Unit}(\tilde{G})$, we have $\phi(F)=\phi(G)$ so by Lemma \ref{connected_again} there is an equivalence $F\to G$ in $\Map^{\adj}(\sk_2(\Adj_{(2,1)}),\CC)$, which we can lift to $\tilde{F}\to\tilde{G}$ in $\Map(\Adj_{(4,1)},\CC)$, because $E_{\Unit}$ is a fibration. Then $\tilde{G}$ is our desired lift of $G$. \end{proof}

Now we prove the fibres of $\phi$ are $1$-connected.

\begin{lemma}\label{1connected_1}

Let $F\in\Map(\sk_3(\Adj_{(3,1)}),\CC)$ and $\alpha:F\to F$ in $\Map(\Unit,\CC)$ such that $\alpha=\Id_F$ on $\{X,Y,l\}$. Then there exists $m:\alpha\to\Id_F$ in $\Map(\{X,Y,l,r\},\CC)$ with $m=\Id^{(2)}_F$ on $\{X,Y,l\}$.
 
\end{lemma}
\begin{proof}Same as \cite[Lemma 3.11]{adj3}\end{proof}

\begin{lemma}\label{1connected_2}

Let $F\in\Map(\sk_4(\Adj_{(4,1)}),\CC)$ and $\alpha:F\to F$ in $\Map(\Unit,\CC)$, such that $\alpha=\Id_F$ on $\{X,Y,l,r\}$. Then there exists $m:\alpha\to\Id_F$ in $\Map(\Unit,\CC)$ with $m=\Id^{(2)}_F$ on $\{X,Y,l\}$.
 
\end{lemma}

\begin{proof}
Same as \cite[Lemma 6.14]{adj3}.\end{proof}

\begin{lemma}\label{1connected_again}
The fibres of $\phi$ are $1$-connected.
\end{lemma}
\begin{proof}
 
Consider $\alpha:F\to F$ in $\Map^{\adj}(\Unit,\CC)$ with $\alpha=\Id_F$ on $\{X,Y,l\}$. We need to construct $m:\alpha\to\Id_F$ in $\Map^{\adj}(\Unit,\CC)$ with $m=\Id^{(2)}_F$ on $\{X,Y,l\}$.

By Lemma \ref{lrusurj} we can extend $F$ to $\Adj_{(4,1)}$.

By Lemma \ref{1connected_1}, we can construct $m_1:\alpha\to\Id_F$ in $\Map(\sk_1(\Adj_{(2,1)}),\CC)$ with $m_1=\Id_F^{(2)}$ on $\{X,Y,l\}$. By Theorem \ref{fibration}, we can extend this to $m_1:\alpha\to\alpha_1$ in $\Map^{\adj}(\Unit,\CC)$, where $\alpha_1=\Id_F$ over $\sk_1(\Adj_{(2,1)})$ and $m_1=\Id^{(2)}_F$ over $\{X,Y,l\}$.

By Lemma \ref{1connected_2} we can construct $m_2:\alpha_1\to\Id_F$ in $\Map(\Unit,\CC)$ with $m_2=\Id^{(2)}_F$ on $\{X,Y,l\}$.

Composing $m_1$ and $m_2$ we obtain the desired $m:\alpha\to\Id_F$ in $\Map^{\adj}(\Unit,\CC)$, with $m=\Id^{(2)}_F$ on $\{X,Y,l\}$. \end{proof}

Now we show the fibres of $\phi$ are $2$-connected.

\begin{lemma}\label{2connected_1}
Let $F\in\Map(\sk_4(\Adj_{(4,1)}),\CC)$ and $m:\Id_F\to\Id_F$ in $\Map(\Unit,\CC)$ such that $m=\Id^{(2)}_F$ on $\{X,Y,l\}$. Then there exists $\A:m\to \Id^{(2)}_F$ in $\Map(\{X,Y,l,r\},\CC)$ with $\A=\Id^{(3)}_F$ on $\{X,Y,l\}$.
\end{lemma}
\begin{proof}This is almost the same as the proof of \cite[Lemma 6.17]{adj3}, except  there we use the $4$-morphism $m_c$ in the construction of $\A_r$. There is a completely analogous construction using $m_u$ instead.\end{proof}

\begin{lemma}\label{2connected_2}

Let $F\in\Map(\Adj_{(4,1)},\CC)$ and $m:\Id_F\to\Id_F$ in $\Map(\Unit,\CC)$ such that $m=\Id^{(2)}_F$ on $\{X,Y,l,r\}$. Then there exists $\A:m\to \Id^{(2)}_F$ in $\Map(\Unit,\CC)$ with $\A=\Id^{(3)}_F$ on $\{X,Y,l\}$. 
 
\end{lemma}

\begin{proof}Same as Lemma \ref{aaudef}.\end{proof}

\begin{lemma}
The fibres of $\phi$ are $2$-connected.
\end{lemma}
\begin{proof}This follows from Lemmas \ref{2connected_1} and \ref{2connected_2}, by an argument which is analogous to the one in the proof of Lemma \ref{1connected_again}.\end{proof}

Now we show the fibres of $\phi$ are $3$-connected.

\begin{lemma}\label{3connected_1}
Let $F\in\Map(\Adj_{(4,1)},\CC)$ and $\A:\Id^{(2)}_F\to\Id^{(2)}_F$ in $\Map(\Unit,\CC)$ such that $\A=\Id^{(3)}_F$ on $\{X,Y,l\}$. Then $\A_r=\id^{(3)}_{F(r)}$.
\end{lemma}
\begin{proof}Same as Lemma \ref{3connected_core}.\end{proof}

\begin{lemma}
The fibres of $\phi$ are $3$-connected.
\end{lemma}
\begin{proof}Follows immediately from Lemma \ref{3connected_1}. \end{proof}

\begin{proposition}
The map $\phi:\Map^{\adj}(\Unit,\CC)\to\Map^L(\theta^{(1)},\CC)$ is a trivial fibration of $4$-groupoids. 
\end{proposition}
\begin{proof}
We know $\phi$ is a fibration of $4$-groupoids. Its fibres are $3$-groupoids because $\theta^{(1)}$ contains the $0$-skeleton of $\Unit$. Since the fibres are $3$-connected, they are weakly contractible, so $\phi$ is a trivial fibration. \end{proof}

\section{A conjectural description of $n$-categorical adjunctions}\label{conjecture}

\begin{definition}

An \textbf{adjoint $k$-equivalence} $X\to Y$ in an $n$-category $\CC$, for $2\leq k\leq n+1$, is an adjoint $1$-equivalence in the $n-k+1$-category $\Hom_{\CC}(s(X),t(X))$. An \textbf{adjoint $1$-equivalence} in an $n$-category is an isomorphism for $n=1$ and an equality for $n=0$. When $n\geq 2$, an adjoint $1$-equivalence $X\to Y$ in an $n$-category $\CC$ consists of the following data.

\begin{enumerate}[start=0,label={(\arabic*)}]
 \item $0$-cells $X=$ \includegraphics[align=c,scale=2]{adj/pres/x.pdf} and $Y=$ \includegraphics[align=c,scale=2]{adj/pres/y.pdf} 
 \item $1$-cells $$f= \includegraphics[align=c,scale=1]{adj/pres/l.pdf} :\xymatrix@1{X\ar@1@<1ex>[r] & \ar@1@<1ex>[l]Y}: \includegraphics[align=c,scale=1]{adj/pres/r.pdf} =f^{-1}$$
 \item two adjoint $2$-equivalences 
 
 \begin{center}
  \begin{tabular}{cccccc}

 $\un=$ & \includegraphics[align=c,scale=1]{adj/pres/u.pdf} & $:$ & \includegraphics[align=c,scale=1]{adj/pres/u_s.pdf} & $\xymatrix@1{\ar@2[r] & }$ & \includegraphics[align=c,scale=1]{adj/pres/u_t.pdf} \\ \\
 
  $\co=$ & \includegraphics[align=c,scale=1]{adj/pres/c.pdf} & $:$ & \includegraphics[align=c,scale=1]{adj/pres/c_s.pdf} & $\xymatrix@1{\ar@2[r] & }$ & \includegraphics[align=c,scale=1]{adj/pres/c_t.pdf} 
 
   \end{tabular}

 \end{center}
 
 \item an adjoint $3$-equivalence
 
  \begin{center}
  \begin{tabular}{cccccc}

 $C_f=$ & \includegraphics[align=c,scale=1]{adj/pres/cusp_l.pdf} & $:$ & \includegraphics[align=c,scale=1]{adj/pres/snake_l.pdf} & $\xymatrix@1{\ar@3[r] & }$ & \includegraphics[align=c,scale=1]{adj/pres/id_l.pdf} 
 
   \end{tabular}

 \end{center}
 
 \end{enumerate}
 
\end{definition}

\begin{conjecture}
 
This notion of adjoint $k$-equivalence is a \textbf{coherent $k$-equivalence}, in the sense that the space of adjoint $k$-equivalences in an $n$-category is equivalent to the space of $k$-morphisms which are equivalences.

\end{conjecture}

\begin{definition}

We define $\Adj_{(n,1)}$ to be the $(n+1)$-computad containing the following cells.

\begin{enumerate}[start=0,label={(\arabic*)}]
 \item $0$-cells $X=$ \includegraphics[align=c,scale=2]{adj/pres/x.pdf} and $Y=$ \includegraphics[align=c,scale=2]{adj/pres/y.pdf} 
 \item $1$-cells $$l= \includegraphics[align=c,scale=1]{adj/pres/l.pdf} :\xymatrix@1{X\ar@1@<1ex>[r] & \ar@1@<1ex>[l]Y}: \includegraphics[align=c,scale=1]{adj/pres/r.pdf} =r$$
 \item $2$-cells 
 
 \begin{center}
  \begin{tabular}{cccccc}

 $\un=$ & \includegraphics[align=c,scale=1]{adj/pres/u.pdf} & $:$ & \includegraphics[align=c,scale=1]{adj/pres/u_s.pdf} & $\xymatrix@1{\ar@2[r] & }$ & \includegraphics[align=c,scale=1]{adj/pres/u_t.pdf} \\ \\
 
  $\co=$ & \includegraphics[align=c,scale=1]{adj/pres/c.pdf} & $:$ & \includegraphics[align=c,scale=1]{adj/pres/c_s.pdf} & $\xymatrix@1{\ar@2[r] & }$ & \includegraphics[align=c,scale=1]{adj/pres/c_t.pdf} 
 
   \end{tabular}

 \end{center}
 
 \item two adjoint $3$-equivalences
 
  \begin{center}
  \begin{longtable}{cccccc}

 $C_l^{-1}=$ & \includegraphics[align=c,scale=1]{adj/pres/cusp_l_inv.pdf} & $:$ & \includegraphics[align=c,scale=1]{adj/pres/id_l.pdf} & $\xymatrix@1{\ar@3[r] &}$ & \includegraphics[align=c,scale=1]{adj/pres/snake_l.pdf} 
 
 \\
 
 \\
 
 $C_r=$ & \includegraphics[align=c,scale=1]{adj/pres/cusp_r.pdf} & $:$ & \includegraphics[align=c,scale=1]{adj/pres/snake_r.pdf} & $\xymatrix@1{\ar@3[r] &}$ & \includegraphics[align=c,scale=1]{adj/pres/id_r.pdf} 
 
   \end{longtable}

 \end{center}
 
 \item an adjoint $4$-equivalence
 
 \begin{tabular}{cccccc}$\SW=$ &\includegraphics[align=c,scale=1]{adj/pres/swallowtail.pdf}  & $:$ & \includegraphics[align=c,scale=0.7]{adj/pres/swallowtail_s.pdf} & $\to$ & \includegraphics[align=c,scale=0.7]{adj/pres/id_u.pdf}\end{tabular}
 
 \end{enumerate}

\end{definition}

\begin{conjecture}
 
Given an $n$-category $\CC$, the space of maps $\Adj_{(n,1)}\to\CC$ is equivalent to the space of $1$-morphisms in $\CC$ which admit a right adjoint. So such a map deserves to be called a \textbf{coherent adjunction}.
 
\end{conjecture}

\end{document}